\documentclass[11pt,a4paper]{amsart}
\date{19 October, 2023}
\usepackage{latexsym,amsmath,amsthm,amsfonts,amscd,amssymb,mathdots,mathtools}
\usepackage[utf8]{inputenc}
\usepackage{stmaryrd}
\usepackage[mathscr]{eucal} 
\usepackage{mathrsfs}
\usepackage{enumerate}
\usepackage[margin=1in]{geometry}
\usepackage{mathrsfs}
\usepackage{amssymb}
\usepackage[all]{xy}
\usepackage{xcolor}
\usepackage{graphics}
\usepackage{lscape}
\usepackage{array}
\usepackage{xy}

\numberwithin{equation}{section}
\newtheorem*{theorem*}{Theorem}
\newtheorem{theorem}{Theorem}[section]
\newtheorem{lemma}[theorem]{Lemma}
\newtheorem{proposition}[theorem]{Proposition}
\newtheorem{corollary}[theorem]{Corollary}
\theoremstyle{definition}

\theoremstyle{remark}
\newtheorem{remark}[theorem]{Remark}

\usepackage{hyperref}\hypersetup{colorlinks}

\usepackage{color} 

\definecolor{darkred}{rgb}{1,0,0} 
\definecolor{darkgreen}{rgb}{0,1,0}
\definecolor{darkblue}{rgb}{0,0,1}

\hypersetup{colorlinks,
linkcolor=darkblue,
filecolor=darkgreen,
urlcolor=darkred,
citecolor=darkgreen}

\DeclareMathOperator{\Jac}{Jac}

\DeclareMathOperator{\Hom}{Hom}

\newcommand{\Nm}{\mathrm{Nm}}

\newcommand{\BBB}{\mathrm{(BBB)}}
\newcommand{\BAA}{\mathrm{(BAA)}}

\newcommand{\Id}{\mathrm{Id}}
\newcommand{\pr}{\mathrm{pr}}
\newcommand{\sing}{\mathrm{sing}}

\newcommand{\supp}{\mathrm{supp}}
\newcommand{\End}{\mathrm{End}}

\newcommand{\Sym}{\mathrm{Sym}}

\newcommand{\GL}{\mathrm{GL}}
\newcommand{\SL}{\mathrm{SL}}
\newcommand{\PSL}{\mathrm{PSL}}

\newcommand{\tr}{\mathrm{tr}}

\renewcommand{\div}{\mathrm{div}}

\newcommand{\Cc}{\mathcal{C}}

\newcommand{\Ee}{\mathcal{E}}

\newcommand{\Ll}{\mathcal{L}}
\newcommand{\Mm}{\mathcal{M}}
\newcommand{\Nn}{\mathcal{N}}
\newcommand{\Oo}{\mathcal{O}}

\newcommand{\Ss}{\mathcal{S}}

\newcommand{\Vv}{\mathcal{V}}
\newcommand{\Ww}{\mathcal{W}}

\newcommand{\M}{\mathbf{M}}

\newcommand{\wt}[1]{\widetilde{#1}}

\newcommand{\CC}{\mathbb{C}}

\newcommand{\ZZ}{\mathbb{Z}}

\newcommand{\HH}{\mathbb{H}}

\newcommand{\PP}{\mathbb{P}}

\renewcommand{\to}{\longrightarrow}

\newcommand{\quotient}[2]{{\raisebox{.2em}{\thinspace $#1$}\left / \raisebox{-.15em}{ $#2$}\right.}}

\newcommand{\smtrx}[1]{\left (\begin{smallmatrix}#1\end{smallmatrix}\right)}

\newcommand\Quotient[2]{
\mathchoice
{
\text{\raise1ex\hbox{\thinspace $#1$}\Big{/} \lower1ex\hbox{$#2$} \thinspace}%
}
{
#1\,/\,#2
}
{
#1\,/\,#2
}
{
#1\,/\,#2
}
}

\newcommand\GIT[2]{
\mathchoice
{
\text{\raise1ex\hbox{\thinspace $#1$}\Big{/}\!\!\!\!\Big{/} \lower1ex\hbox{$#2$} \thinspace}%
}
{
#1\,/\,#2
}
{
#1\,/\,#2
}
{
#1\,/\,#2
a       }
}

\newcommand{\map}[5]{\begin{array}{ccc}   #1  & \stackrel{#5}{\longrightarrow} &  #2  \\  #3 & \longmapsto & #4  \end{array}}

\title[Narasimhan--Ramanan branes and wobbly Higgs bundles]{\bf Narasimhan--Ramanan branes and wobbly Higgs bundles}

\author[E. Franco]{Emilio Franco}
\address{E. Franco,
\newline\indent Depto. Matem\'aticas, Facultad de Ciencias, 
\newline\indent Universidad Aut\'onoma de Madrid
\newline\indent Campus de Cantoblanco 28049, Madrid, Espa\~na.}
\email{emilio.franco@uam.es}

\author[P.~B.\ Gothen]{Peter B.\ Gothen}
\address{P.~B. Gothen,
\newline\indent Centro de Matemática da Universidade do Porto, 
\newline\indent Departamento de Matemática, Faculdade de Ci\^encias, Universidade do Porto, 
\newline\indent Rua do Campo Alegre s/n, 4169-007 Porto, Portugal}
\email{pbgothen@fc.up.pt}

\author[A. Oliveira]{Andr\'e Oliveira}
\address{A. Oliveira, 
\newline\indent Centro de Matem\'atica da Universidade do Porto,
\newline\indent Departamento de Matemática, Faculdade de Ci\^encias, Universidade do Porto, 
\newline\indent
Rua do Campo Alegre s/n, 4169-007 Porto, Portugal \newline\indent and\newline\indent Departamento de Matem\'atica, Universidade de Tr\'as-os-Montes e Alto Douro, UTAD,\newline\indent
Quinta dos Prados, 5000-911 Vila Real, Portugal}
\email{andre.oliveira@fc.up.pt\newline\indent agoliv@utad.pt}

\author[A. Pe\'on-Nieto]{Ana Pe\'on-Nieto}
\address{A. Pe\'on-Nieto,
\newline\indent
Department of Mathematics, Universidade de Santiago de Compostela, Spain
\newline\indent 
and\newline\indent
	School of Mathematics, University of Birmingham,  UK, 
}
\email{ana.peon@usc.es, a.peon-nieto@bham.ac.uk}

\thanks{%
First author supported by the Scientific Employment Stimulus program, funded by FCT (Portugal), by an FCT Investigator grant with fellowship reference CEECIND/04153/2017. Second and third authors partially supported by CMUP (UIDB/00144/2020 and UIDP/00144/2020) and the project EXPL/MAT-PUR/1162/2021 funded by FCT (Portugal) with national funds. Fourth author partially supported by the  project ReaCH, funded through the Consolidaci\'on Investigadora Programme, grant number CNS2022-136042, the Marie Sklodowska-Curie European Individual Fellowship GoH, grant number 897722, and the COVID Support Programme from the University of Birmingham.
}

\begin{document}

\begin{abstract}
Narasimhan--Ramanan branes, introduced by the authors in a previous article, consist of a family of $\BBB$-branes inside the moduli space of Higgs bundles, and a family of complex Lagrangian subvarieties. It was conjectured that these complex Lagrangian subvarieties support the $\BAA$-branes that are mirror dual to the Narasimhan--Ramanan $\BBB$-branes. In this article we show that the support of these branes intersects non-trivially the locus of wobbly Higgs bundles. 
\end{abstract}

\maketitle

\begin{flushright}
{\small{\it To Óscar García-Prada,\\
on the occasion of his 60th anniversary.}} \\ \vspace*{0.25cm}
\end{flushright}

\section{Introduction}

Let $\Mm$ denote the moduli space of rank $r$ and degree $d$ Higgs bundles over a smooth complex projective curve $X$. Among the many important contributions of Óscar García-Prada's is his work with S.~Ramanan on the description of finite order automorphisms of the moduli space $\Mm$ and their fixed-point loci; see \cite{garciaprada-ramanan-2,garciaprada-ramanan}. An important instance of such an automorphism is obtained by tensoring a Higgs bundle by a fixed line bundle $T\to X$ of order $r$. 
The version of this automorphism on the moduli space of vector bundles (without Higgs field) was first considered by M.S. Narasimhan and S. Ramanan in \cite{NR}. They gave a beautiful description of the fixed point locus in the moduli space of bundles in terms of line bundles on the $r$-sheeted unramified cover of the curve $X$ corresponding to the order $r$ line bundle $T$ (see also \cite{garciaprada-ramanan-2,nasser:2005}). A similar picture arises for the fixed point locus $\Mm^T$ in the moduli space of Higgs bundles, as first observed by T. Hausel and M. Thaddeus \cite{HT} in their proof of topological mirror symmetry between the $\SL$- and the $\PSL$-Higgs bundle moduli spaces, where the subspaces $\Mm^T$ play a crucial role.

With this in mind, in \cite{FGOP} we study the fixed point subvarieties $\Mm^T\subset\Mm$ from the point of view of mirror symmetry. The moduli space $\Mm$ is hyperkähler and we prove that the subvarieties $\Mm^T$ are hyperkähler subvarieties of $\Mm$ which support hyperholomorphic line bundles. Such a gadget --- a hyperholomorphic line bundle over a hyperkähler subvariety --- is usually called a $\BBB$-brane on $\Mm$. Now, roughly speaking, mirror symmetry predicts the existence of a dual object, a so called $\BAA$-brane, which consists of flat bundle over a complex Lagrangian subvariety of $\Mm$. Moreover, these dual branes should be related by a Fourier-Mukai transform on $\Mm$ relative to the Hitchin map $h:\Mm\to B$. In \cite{FGOP} we proved that if $d=0$ (so that $h$ admits a section) then indeed the $\BBB$-branes supported on $\Mm^T$ are transformed under Fourier-Mukai into a sheaf supported on a suitable Lagrangian of $\Mm$. As predicted by mirror symmetry, these Lagrangians depend on the hyperholomorphic bundle over $\Mm^T$ and have a description in terms of spectral data of the corresponding Higgs bundles. From the point of view of objects on $X$, the corresponding Higgs bundles can also be conveniently described in terms of Hecke transformations.

In \cite{FGOP} the Fourier-Mukai transform is taken fiberwise with respect to $h$ for a generic fiber of $h$ over the image $h(\Mm^T)$, and so the description of the Lagrangians in \cite{FGOP} is not complete. In particular their intersection with the nilpotent cone $h^{-1}(0)$ is not described there. 
Here we aim to give a first step in that direction, by computing the limits at zero of the $\CC^*$-flows (for the standard $\CC^*$-action on $\Mm$ scaling the Higgs field) in those Lagrangians when the rank $r$ equals $2$. These limits lie in $h^{-1}(0)$ and are $\CC^*$-fixed points. In the language of T. Hausel and N. Hitchin \cite{HH} a stable Higgs bundle  is called wobbly if it is the limit of a non-isomorphic nilpotent Higgs bundle, and very stable otherwise.

On the other hand the intersection of the $\BBB$-branes with the nilpotent cone is known and lies in the moduli space of bundles (Higgs bundles with zero Higgs field). Here we analyse their intersection with the components of the wobbly locus.

We show that the limits of Higgs bundles in our Lagrangians are in fact wobbly Higgs bundles (with the obvious exception of limits which are fixed points lying at the top of the nilpotent cone). In particular, our Lagrangians are of a complementary nature to those considered in \cite{HH}, which arise as upward flows from very stable Higgs bundles.

This paper is organized as follows. In Section \ref{sec:moduli-spaces} we provide a short introduction to Higgs bundles and their moduli spaces, we recall the mirror symmetry conjecture for branes and review the definition of very stable and wobbly bundles and Higgs bundles. Section \ref{sec:NR-branes} is a survey of the results obtained in \cite{FGOP} focusing on the special case of rank two bundles. Hence, in this section we recall the construction of the Narasimhan--Ramanan $\BBB$-branes, their associated spectral data and the complex Lagrangian subvarieties that conjecturally support the dual Narasimhan--Ramanan $\BAA$-branes. Finally, we present original results in Section \ref{sc NR and wobbly}. We first describe, in Section \ref{sc BBB and wobbly}, how the support of the Narasimhan--Ramanan $\BBB$-branes hits the locus of wobbly vector bundles (cf.~Theorem \ref{thm:wobbly_BBB}). Secondly, in Section \ref{sc stratification}, we construct a certain stratification of the generic intersection of our complex Lagrangian subvarieties with a Hitchin fibre that will be used, in Section \ref{sc C*flows of lagrangians}, to compute their $\mathbb{C}^*$-fixed points. It follows from this analysis that most of the times these complex Lagrangians intersect the $\CC^*$-fixed points loci only in wobbly Higgs bundles, the exception being the intersection with the $\CC^*$-fixed point loci corresponding to maximal degree (the ``tip of the nilpotent cone'') where the wobbly locus is empty (see Theorem \ref{thm:wobblylimits} for the detailed statement).

\subsection*{Acknowledgement.} We thank Johannes Horn for useful discussions.

\section{Moduli spaces of Higgs bundles}
\label{sec:moduli-spaces}

\subsection{Higgs bundles and their moduli spaces}

Let $X$ be a smooth projective complex curve of genus $g>1$. A Higgs bundle of rank $r$ and degree $d$ over $X$ is a pair $(E,\varphi)$ where $E$ is a rank $r$ and degree $d$ holomorphic vector bundle over $X$ and the Higgs field $\varphi$ is a holomorphic map $\varphi:E\to E\otimes K$, with $K$ the canonical line bundle of $X$.

Recall that a vector bundle $E$ is semistable (resp.~stable) if the slope of any non-zero proper subbundle is less than or equal to (resp.~strictly less) to the slope of $E$. The (semi)stability condition for a Higgs bundle $(E,\varphi)$ is the same as for $E$, but only applied to subbundles which are $\varphi$-invariant. In both cases an object which is not semistable is said to be unstable. 
Note that $(E,\varphi)$ may be semistable even if its underlying vector bundle $E$ is unstable.

Let  $\Mm(r,d)$ (resp.~$\Nn(r,d)$) be the moduli space of semistable Higgs (resp.~vector) bundles over $X$ of rank $r$ and degree $d$. The locus $\Mm(r,d)^s$ (resp.~$\Nn(r,d)^s$) of stable Higgs (resp.~vector) bundles is smooth. Considering a vector bundle as a Higgs bundle with zero Higgs field we have an inclusion $\Nn(r,d)\subset\Mm(r,d)$. Moreover, since the cotangent space $H^1(X,\End(E))^*$ at $E\in\Nn(r,d)^s$ is isomorphic to the space $H^0(X,\End(E)\otimes K)$ of Higgs fields on $E$ by Serre duality, it follows that the cotangent bundle $T^*\Nn(r,d)^s$ of the stable locus is an open dense subspace of $\Mm(r,d)$.

For the remainder of the paper our main focus will be on the case of rank $2$ and degree $0$ (even though everything we say in Sections~\ref{sec:moduli-spaces} and \ref{sec:NR-branes}  can be done for arbitrary rank). So let now \[\Mm=\Mm(2,0),\ \ \Nn=\Nn(2,0).\] We will need at some point to consider rank $1$ Higgs bundles of degree $0$. Any rank 1 Higgs bundle is of course stable and in this case the inclusion of the cotangent bundle of the Jacobian is an isomorphism: $\Mm(1,0)=T^*\Jac(X)=\Jac(X)\times H^0(X,K)$.

Over the smooth locus $T^*\Nn^s$ carries a canonical holomorphic symplectic form $\Omega$ which extends in a natural way to the smooth locus of $\Mm$, which is thus a holomorphic symplectic manifold. Then $\Omega$ gives rise to a holomorphic volume form on $\Mm$, showing that $\Mm$ is a (singular) non-compact Calabi-Yau.

One of the most important features of $\Mm$ is its structure as an integrable system. This comes from the Hitchin map 
\[h:\Mm\longrightarrow B=H^0(X,K)\oplus H^0(X,K^2), \ \ \ \ h(E,\varphi)=(\tr(\varphi),\det(\varphi))\]
which is surjective and proper.
We have that $2\dim(B)=\dim(\Mm)$, and for each $c=(a,b)\in B$ there is a corresponding spectral curve $X_c$ lying in the surface $|K|$ given by the total space of $K$. If $c=h(E,\varphi)$ then $X_c$ is given by the eigenvalues of $\varphi$. Precisely, if $\pi:|K|\to X$ is the canonical projection and $\lambda$ the tautological section of $\pi^*K$, then $X_c$ is the zero locus of the section $\lambda^2+\pi^*a\lambda+\pi^*b\in H^0(|K|,\pi^*K^2)$. From its very definition, one sees that any spectral curve $X_c$ lie in the linear system of $2X=X_0$ and that $\pi_c:X_c\to X$ is a ramified $2$-cover. For generic $c\in B$, the spectral curve $X_c$ is smooth but it can be singular, and even reducible or non-reduced (like $X_0$). The locus of elements $c\in B$ for which $X_c$ is not smooth is called the discriminant locus. 
For $X_c$ integral, one can define its compactified Jacobian $\overline{\Jac}^{2g-2}(X_c)$ as the moduli space of rank $1$ torsion-free sheaves of degree $2g-2$ over $X_c$ (which by irreducibility of $X_c$ are always stable). The compactified Jacobian contains $\Jac^{2g-2}(X_c)$ as a dense open subset parameterizing degree $2g-2$ line bundles on $X_c$. If $X_c$ is smooth, then $\overline{\Jac}^{2g-2}(X_c)=\Jac^{2g-2}(X_c)$. As proved in \cite{BNR,hitchin_duke}, each $\mathcal{F}\in\overline{\Jac}^{2g-2}(X_c)$ determines a Higgs bundle in $h^{-1}(c)$ by taking $(\pi_*\mathcal{F},\pi_*\lambda)$. If $(E,\varphi)=(\pi_*\mathcal{F},\pi_*\lambda)\in h^{-1}(c)$, the spectral curve $X_c$ encodes eigenvalues of $\varphi$ whereas the sheaf $\mathcal{F}$ encodes the corresponding eigenspaces. Hence $(X_c,\mathcal{F})$ is called the spectral data of $(E,\varphi)$. The spectral correspondence gives an isomorphism $\overline{\Jac}^{2g-2}(X_c)\cong h^{-1}(c)$ and hence the generic fiber of $h$ is a torsor for an abelian variety. Moreover, being the fiber of an integrable system, it is in fact Lagrangian in $\Mm$.

As proved in \cite{hitchin-self} the $C^\infty$ manifold $\M$ underlying (the smooth locus of) $\Mm$ is hyperkähler. Thus it has complex structures $I,J$ and $K=IJ$ 
verifying the quaternionic relations and a Riemannian metric which is Kähler with respect to each of them. By convention $I$ is the complex structure of $\Mm$, i.e.~$\Mm=(\M,I)$. Then $(\M,J)$ is the moduli space of flat connections on the $C^\infty$-trivial rank 2 vector bundle over $X$, i.e., the moduli of rank $2$ local systems. Denote it by $\Mm_{dR}=(\M,J)$, the index standing for `de Rham moduli space', a term coined by C. Simpson (who also called $\Mm$ the Dolbeault moduli space) \cite{Si94}. In fact $\M$ has a $2$-sphere of complex structures $aI+bJ+cK$ with $a^2+b^2+c^2=1$ and, with the exception of $\pm I$, all these complex structures are equivalent to $J$.  Write $\omega_I=g(I,-)$, $\omega_J=g(I,-)$ and $\omega_K=g(K,-)$ for the corresponding Kähler forms; then the holomorphic symplectic form $\Omega$ on $\Mm$ mentioned above is given by $\Omega=\omega_J+i\omega_K$.

\subsection{Branes in moduli spaces of Higgs bundles and mirror symmetry}

In this section we give a brief overview of a few important mathematical ideas and concepts coming from mirror symmetry. Our treatment is intended to motivate what follows and is neither complete nor rigorous.

An A-brane on a Kähler manifold is an object of its Fukaya category; an example of an A-brane is a Lagrangian with a flat bundle over it. On the other hand, a B-brane is an object of its derived category of coherent sheaves, for instance a holomorphic vector bundle over a complex submanifold. Since $\M$ is hyperkähler, we have branes for each complex and symplectic structure. Thus a $\BAA$-brane on $\M$ is at the same time a B-brane for $I$ and an A-brane for $\omega_J$ and $\omega_K$. For instance, a flat bundle supported on a complex Lagrangian of $\Mm$ (with respect to $\Omega)$ is such a gadget. In turn, a $\BBB$-brane is a B-brane for $I$, $J$ and $K$, for example a hyperholomorphic bundle over a hyperkähler subvariety of $\Mm$.

Higgs bundles can be considered for any reductive group $G$ and 
R. Donagi and T. Pantev \cite{DP} (motivated by the previous works \cite{HT,hitchin_G2}) showed that the semi-classical limit of mirror symmetry holds for pairs of Higgs bundle moduli spaces for Langlands dual groups in the complement of the discriminant locus. This means that there is a duality of fibers of the Hitchin maps for $G$ and its Langlands dual, realized by a Fourier-Mukai transform.
Now, the moduli space $\Mm$ is the moduli space of $G$-Higgs bundles for $G=\GL(2,\CC)$. Since this group is self-dual (or, equivalently, since the generic fiber of $h$ is a torsor for the self-dual abelian variety $\Jac(X_c)$), then the semi-classical limit of mirror symmetry holds for $\Mm$ in the complement of the discriminant locus.

According to the homological mirror symmetry conjecture proposed by M. Kontsevich in \cite{kont:1995} there should be an equivalence between the Fukaya category of $\Mm_{dR}$ and its derived category of coherent sheaves, and mirror symmetry should be mathematically realized by this equivalence. This equivalence of categories should exchange, in particular, A-branes and B-branes.

Now, C. Simpson proved in \cite{Si94} that $\Mm$ and $\Mm_{dR}$ are, respectively, the fiber over $0$ and $1$ of a holomorphic family of moduli spaces over $\CC$ (called the Hodge moduli space) such that any fiber other than $\Mm$ is isomorphic to $\Mm_{dR}$. Hence mirror symmetry holds for any fiber over $\lambda\in\CC^*$. The so-called semi-classical limit is then $\Mm$, the fiber over $0$ and, as mentioned above, here mirror symmetry amounts to a fiberwise Fourier-Mukai transform relative to the Hitchin map $h$ away from the discriminant locus. Being an algebraic map on $\Mm$, this Fourier-Mukai transform preserves B-branes on $\Mm$. But according to A. Kapustin and E. Witten \cite{kapustin&witten}, after a hyperkähler rotation, it should exchange a B-brane for complex structure $J$ (resp.~$K$) with an A-brane for complex structure $K$ (resp.~$J$). Summing up, mirror symmetry in $\Mm$, outside the discriminant locus, is supposed to be realized by a fiberwise Fourier-Mukai transform and it should interchange $\BBB$-branes and $\BAA$-branes whose supports share the same image under $h$. 

The behavior of mirror symmetry on branes mapping under $h$ to the discriminant locus is far from being understood. The first example of a pair $\BBB$- and $\BAA$-branes on $\M$ lying over the discriminant locus was studied in \cite{FP}. These branes are actually supported over the locus whose corresponding spectral curves are reducible. Hence, the compactified Jacobians are just coarse moduli spaces rather than fine. For this reason a complete description of a Fourier-Mukai transform is not available, even though in \cite{FP} evidence was given that the pair of branes under consideration should be related under a Fourier-Mukai transform. In \cite{FGOP} (partly motivated by \cite{FP}) we considered another type of branes, still supported over the discriminant locus, but such that the spectral curves are generically irreducible and reduced, i.e.~integral. For such curves the compactified Jacobian is fine, and D. Arinkin proved \cite{arinkin} the existence of a Poincaré sheaf making it possible to perform a Fourier-Mukai transform. This allowed us to explore further the duality the pairs of branes, and hence to provide insight on the classical limit of mirror symmetry works over the discriminant locus. Section \ref{sec:NR-branes} will detail the most important aspects of \cite{FGOP} in the rank 2 case.

\subsection{Very stable and wobbly Higgs bundles}
\label{sc very stable}

By definition the nilpotent cone in $\Mm$ is the preimage $h^{-1}(0)\subseteq \Mm$ of zero under the Hitchin map. Equivalently, the nilpotent cone is the locus of Higgs bundles with nilpotent Higgs field.
Consider the $\CC^*$-action on $\Mm$ scaling the Higgs field: $t\cdot(E,\varphi)=(E,t\varphi)$, $t\in\CC^*$. It is well-known \cite{hitchin-self, simpson-cstar} that the  $\CC^*$-fixed points in $\Mm$ lie in the nilpotent cone and are either semistable vector bundles i.e.~in $\Nn\subset h^{-1}(0)$ or Higgs bundles of the form 
\begin{equation}\label{eq:nonzerofixedpoint}
\left(L_1\oplus L_2,\begin{psmallmatrix}0 & 0\\ \phi & 0
\end{psmallmatrix}\right),\ \text{ with }\ 1\leq \deg(L_1)\leq g-1\ \text{ and }\ \phi\in H^0(L_1^{-1}L_2K)\setminus\{0\}.
\end{equation}
Moreover, for any semistable Higgs bundle $(E,\varphi)$ in $\Mm$, the limit $\lim_{t\to0}(E,t\varphi)$ exists \cite{Si94} in $\Mm$ and is a $\CC^*$-fixed point. It is $E\in\Nn\subset\Mm$ precisely when $E$ is itself semistable; otherwise is of the form \eqref{eq:nonzerofixedpoint} with $L_1\subset E$ the maximal destabilizing subbundle, $L_2$ the quotient $E/L_1$ and $\phi$ the composition $L_1\hookrightarrow E\xrightarrow{\ \varphi\ }E\otimes K\to L_2\otimes K$.   

We will also make use of the following notions. A semistable vector bundle $E\in\Nn$ is said to be very stable if there are no non-zero nilpotent Higgs fields $E\to E\otimes K$. In other words, $E=(E,0)$ is not the limit of a $\CC^*$-flow contained in $h^{-1}(0)$. If $E$ is not very stable, it is said to be wobbly. Let $E$ be wobbly, and let $\varphi$ be a non-zero nilpotent Higgs field. Then $E$ can be written as an extension of line bundles
\begin{equation} \label{eq filtration wobbly}
0 \longrightarrow F_0 \stackrel{i}{\longrightarrow} E \stackrel{j}{\longrightarrow} F_1 \longrightarrow 0,
\end{equation}
where $F_0=\ker(\varphi)$ and $F_1=E/F_0$. Then $ F_1 = F_0K(-D)$ for some effective divisor $D$ of degree $\delta$. Indeed, $\varphi$ factors as 
\begin{equation} \label{eq wobbly Higgs field}
\varphi = (i\otimes\Id_K) \circ \beta \circ j : E \longrightarrow E \otimes K.
\end{equation}
for some non-zero $\beta \in H^0(F_1^*F_0K)$ whose divisor is $D$. 
Note also that 
\begin{equation} \label{eq det E wobbly}
\det(E) \cong F_0^2K(-D).
\end{equation}

Suppose now that $\deg(E) = 0$ and $E$ is semistable. If $\deg F_0=d_0$, it follows from \eqref{eq det E wobbly} that 
\begin{eqnarray*}
\delta \equiv 0 \! \! \mod 2,
\\
d_0 = 1-g + \delta/2,
\\
0 \leq \delta \leq 2g-2.
\end{eqnarray*}
Let $\Ww_\delta\subset\Nn$ denote the locus of rank semistable vector bundles $E$ fitting in an exact sequence of the form \eqref{eq filtration wobbly} with $\deg(D) = \delta$ and $\deg(F_0) = 1 - g - \delta/2$. The wobbly locus $\Ww\subset\Nn$ decomposes \cite{PalPauly} into 
\[
\Ww = \Ww_0 \cup \dots \cup \Ww_{2g-2}.
\]
One can prove that $\Ww_\delta$ is birational to the projective bundle over $\Sym^\delta(X)\times\Jac^{g-1 - \delta/2}(X)$ with fibre over $(D,F_0)$ equal to $\PP(H^1(\Oo(D)K^*))$. The points in this bundle are then in $1:1$ correspondence with extensions of the form \eqref{eq filtration wobbly}.

Denote as well $\delta_{max}=g-1$ if $g$ is odd and $\delta_{max}=g-2$ otherwise. It was proven in  \cite{PalPauly} that the components associated to $\delta_{max} < \delta \leq 2g-2$ are embedded in $\Ww_{\delta_{max}}$, 
\[
\Ww_\delta \subset \Ww_{\delta_{max}}, \qquad \delta > \delta_{max},
\]
while for $0 < \delta \leq \delta_{max}$ the corresponding components $\Ww_\delta$ have codimension $1$ in the moduli space of vector bundles $\Nn$ and intersect in higher codimension. Hence the wobbly locus is the union of divisors
\[
\Ww = \Ww_0 \cup \dots \cup \Ww_{\delta_{max}}.
\]

In \cite{pauly&peon-nieto} it was proved that $E$ is very stable if and only if the subspace of all Higgs bundles in $\Mm$ flowing to $E$ --- that is the space of Higgs bundles with underlying bundle $E$ --- is closed in $\Mm$.

These notions were generalized to all $\CC^*$-fixed points in $\Mm$ in \cite{HH} as follows. Let $(E,\varphi)\in h^{-1}(0)$ be a Higgs bundle which is fixed under $\CC^*$, so either $\varphi=0$ or is of the form \eqref{eq:nonzerofixedpoint}. Consider the subspace $W^+_{(E,\varphi)}=\{(E',\varphi')\in\Mm\,|\,\lim_{t\to 0}(E',t\varphi')=(E,\varphi)\}$, usually called the upward flow from $(E,\varphi)$. Then $(E,\varphi)$ is said to be very stable if  $(E,\varphi)$ is the only Higgs bundle in $W^+_{(E,\varphi)}$ with nilpotent Higgs field. Otherwise $(E,\varphi)$ is again said to be wobbly. Like in the bundle case, $(E,\varphi)$ is very stable if and only if $W^+_{(E,\varphi)}$ is closed in $\Mm$. 
\begin{remark}\label{rmk:nomaximalwobbly}
    Note that if $(E,\varphi)=(
    L_1\oplus L_2,\varphi)$ is a fixed point of the form \eqref{eq:nonzerofixedpoint} and wobbly then it must be the limit $\lim_{t\to 0}(E',t\varphi')$  of a Higgs bundle $(E',\varphi')\in h^{-1}(0)$ such that $\lim_{t\to\infty}(E',t\varphi')=(
    L_1'\oplus L_2',\psi)$ is still of the form \eqref{eq:nonzerofixedpoint} but with kernel $L_1'$ such thar $\deg(L_1')>\deg(L_1)$. This means that there are no wobbly Higgs bundles of the form \eqref{eq:nonzerofixedpoint} when $\deg(L_1)$ is maximal i.e.~$\deg(L_1)=g-1$. In all other cases, the wobbly locus is non-empty.
\end{remark}

Upwards flows $W^+_{(E,\varphi)}$ are Lagrangian in $\Mm$, hence with their structure sheaf become $\BAA$-branes which were studied in \cite{HH} in the context of mirror symmetry. In section 8.2 of loc.~cit.~the authors briefly consider upward flows from wobbly Higgs bundles. We will see below an example of Lagrangians in $\Mm$ whose corresponding points always flow onto wobbly Higgs bundles except when such limits are of the form \eqref{eq:nonzerofixedpoint} with $\deg(L_1)=g-1$ (which we have seen must be very stable).

\section{Narasimhan--Ramanan branes}
\label{sec:NR-branes}

\subsection{$\BBB$-branes out of fixed point loci under tensorization}

Consider the group
$\Jac(X)[2]=\{T\in \Jac(X)\,|\, T^2\cong\Oo_X\}\cong (\ZZ/2)^{2g}$ of $2$-torsion line bundles on $X$, acting on $\Mm$ by tensorization $(E,\varphi)\cdot T=(E\otimes T,\varphi)$. 

Fix, once and for all, a non-trivial $T\in\Jac(X)[2]\setminus\{\Oo_X\}$ and consider the fixed point set under tensorization by $T$,
\[\Mm^T=\{(E,\varphi)\in\Mm\,|\,(E,\varphi)\cong(E\otimes T,\varphi)\}\subset\Mm.\] It is clearly a complex subvariety of $\Mm=(\M,I)$, but since $T$ is flat and tensorization is holomorphic also in complex structure $J$, $\Mm^T$ is actually a hyperkähler subvariety of $\Mm$. Moreover, following the arguments of \cite{NR} (see also \cite{nasser:2005,garciaprada-ramanan-2,HT}) one can check that $\Mm^T$ is explicitly characterized using rank $1$ Higgs bundles over the étale $2$-cover $p:X_T\to X$ canonically defined by the line bundle $T$. Namely, pushforward under $p$ induces an isomorphism 
\begin{equation}\label{eq:isom}
\Mm^T\cong \quotient{T^*\Jac(X_T)}{\ZZ_2},
\end{equation}
where $\ZZ_2$ denotes the Galois group of $p$ acting by pullback; see \cite[Theorem 3.19]{FGOP} and also \cite[Proposition 7.1]{HT}.

\begin{remark}
Fix a degree $m$ non-trivial étale cover $p:C\to X$. In \cite{FGOP} we treated the more general case of hyperkähler subvarieties of the rank $r$ and degree $0$ Higgs bundles moduli space which parameterize semistable Higgs bundles arising as pushforward under $p$ of degree $0$ semistable Higgs bundles of rank $s$ over $C$, where $sm=r$. According to \eqref{eq:isom}, $\Mm^T$ is an example of these subvarieties with $C=X_T$, $s=1$ and $m=2$. The case of trivial coverings is the subject of \cite{FP}.
\end{remark}

Now, out of any degree $0$ line bundle $\Ll\in\Jac(X)$, we build a hyperholomorphic line bundle on $\Mm^T$, using \eqref{eq:isom}, as follows. Let $\Delta:\Mm^T\to\Jac(X)$ be defined as 
$\Delta(E,\varphi)=\Nm_p(F)$ where $F\in\Jac(X_T)$ is such that $E=p_*F$ and $\Nm_p$ is the norm map (which is invariant under the Galois group).
From the self-duality of the Jacobian of $X$, $\Ll$ uniquely determines a flat line bundle 
 $\check{\Ll}\to \Jac(X)$. We then take the line bundle $\Delta^*\check\Ll$ over $\Mm^T$, which is hyperholomorphic because it is flat (hence its curvature is trivially of type $(1,1)$ with respect to both $I$ and $J$).  Hence, we have the following
 
\begin{theorem}\label{thm:BBBbrane}
The line bundle $\Delta^*\check\Ll$ is a $\BBB$-brane of $\Mm$ supported on $\Mm^T$.
\end{theorem} 

In \cite{FGOP} we called this a Narasimhan-Ramanan $\BBB$-brane on $\Mm$ due to their pioneering work \cite{NR}. We are interested in knowing how $\Delta^*\check\Ll$ behaves under mirror symmetry. As briefly explained above the classical limit of mirror symmetry in $\Mm$ is realized by a relative Fourier-Mukai under $h$, so the idea is to look at the spectral data of the Higgs bundles in $\Mm^T$ and then perform a Fourier-Mukai transform of the sheaf $\Delta^*\check\Ll\cap h^{-1}(c)$ for $c$ in the image of $h$ restricted to $\Mm^T$. 

\subsection{The corresponding spectral curves}

Let $B^T=h(\Mm^T)$ be the image of $\Mm^T$ under $h$; it is a closed subset of $B$. In order to obtain an explicit description of $B^T$, consider 
\[
\wt B^T := H^0(X,K) \oplus H^0(X,KT),
\]
which is naturally isomorphic to $H^0(X_T,K_T)$ thanks to the pushforward under $p : X_T \to X$. As stated in \cite[(3.14)]{FGOP} (see also \eqref{eq:isom}), $B^T$ amounts to the quotient
\[
B^T = \quotient{\wt B^T}{\ZZ_2},
\]
where the Galois group $\ZZ_2$ acts trivially on $H^0(X, K)$ and with negative sign on $H^0(X,TK)$. The embedding of $B^T$ into $B$ can be described as follows, first consider the morphism
\[
\map{\wt B^T = H^0(X,K) \oplus H^0(X,KT)}{B = H^0(X,K) \oplus H^0(X, K^2)}{(a, b)}{c = (2a, a^2 - b^2),}{}
\]
noting that two points in $\wt B^T$ have the same image if and only if they are related by the action of the Galois group $\ZZ_2$. Then, the above morphism factors through
\[
\map{B^T}{B}{(a, \pm b)}{c = (2a, a^2 - b^2),}{}
\]
which is injective. Since $T \cong T^*$, observe that $b \in H^0(X,TK)$ can be understood as a map between the total spaces  
\[
b : |T| \to |K|.
\]
The smooth curve $X_T$ arises as a spectral curve, in the total space $|T|$, associated to the non-vanishing section $1$ of $T^2 \cong \Oo_X$. It then follows that $X_c$ can be described as the image of $X_T$ under $b$ shifted by $a$,
\begin{equation} \label{eq X_c from b}
X_c = \nu_{(a,b)}(X_T),
\end{equation}
where $\nu_{(a,b)}(\bullet)$ denotes $b(\bullet) + a$. Observe that $\nu_{(a,b)}$ is invariant with respect to the Galois $\ZZ_2$-action. Since $X_T$ is invariant under the Galois group $\ZZ_2$, one has that $\nu_{(a,b)}$ and $\nu_{(a,-b)}$ have the same image $X_c$.

It then follows that $X_c$ is irreducible and singular so no smooth spectral curves are parameterized by $B^T$. If $X_c$ is integral and given by $c=(2a,a^2-b^2)$ then it is easy to see that its singularities project onto the double zeros of $b^2$ hence to the zeros of $b$, thus \[
\pi_c\left ( \sing(X_c) \right ) = \div(b),
\] with $\div(b)$ the divisor of $b$, so that $\Oo(\div(b))\cong KT$. Notice that $\div(b)$ is a simple divisor because the the singularites of $X_c$ are ordinary double points.

Let $B^T_{ni}\subset B^T$ denote the open dense subset parameterizing spectral curves in $B^T$ which are nodal and integral. For $c=(2a,a^2-b^2)\in B_{ni}^T$, $X_c$ has exactly $2g-2$ ordinary double points as singularities and $\pi_c:X_c\to X$ is ramified precisely at those points. Furthermore, since $X_T$ is smooth and $\nu_{(a,b)}:X_T\to X_c$ surjects onto $X_c$, being a generic isomorphism (namely, outside the points projecting to $\div(b)$), it is the normalization morphism. 

All of the above is summarized in the following theorem.

\begin{theorem}{\cite[Theorem 3.11]{FGOP}}\label{thm:spectral curves}
Let $c=(2a,a^2-b^2)\in B^T_{ni}$ and $X_c$ the corresponding nodal spectral curve with $\pi_c:X_c\to X$ the corresponding ramified $2$-cover. Then $X_c$ is normalized by $X_T$ and if $\nu_{(a,b)}:X_T\to X_c$ is the normalization morphism, we have $p=\pi_c\circ\nu_{(a,b)}$. Moreover, the nodes of $X_c$ project onto the simple effective divisor of $b\in H^0(X,KT)$.
\end{theorem}

We have an embedding 
\[\check\nu_{(a,b)}:\Jac(X_T)\hookrightarrow\overline{\Jac}^{2g-2}(X_c)\] given by pushforward under $\nu_{(a,b)}$ (whose image is contained in $\overline{\Jac}^{2g-2}(X_c)\setminus\Jac^{2g-2}(X_c)$). From this we have the spectral data parameterizing Higgs bundles in $\Mm^T$ and also the description of the restriction to $\Mm^T\cap h^{-1}(c)$, for each $c\in B^T$, of the $\BBB$-brane constructed in Theorem \ref{thm:BBBbrane}.

\begin{theorem}{\cite[Propositions 3.13 and 4.6]{FGOP}}\label{thm:spectral data fixedpoints}
For every $c\in B^T_{ni}$, we have $\Mm^T\cap h^{-1}(c)=\mathrm{Im}(\check\nu_{(a,b)})\cong\Jac(X_T)$. Furthermore, the restriction $\Delta^*\check\Ll|_{\Mm^T\cap h^{-1}(c)}$ is isomorphic to the pushforward of $\Nm_p^*\check\Ll\in\Jac(X_T)$ under $\check\nu_{(a,b)}$.
\end{theorem}

Both Theorems \ref{thm:spectral curves} and \ref{thm:spectral data fixedpoints} follow by analyzing the spectral data of Higgs bundles fixed under $T$, using the isomorphism \eqref{eq:isom}.

The restriction $\Delta^*\check\Ll|_{\Mm^T\cap h^{-1}(c)}$ is then the bundle we want to Fourier-Mukai, but let us before introduce another subvariety of $\Mm$, constructed in terms of the spectral correspondence.

\subsection{The Lagrangians}

Fix $c\in B^T_{ni}$ and consider then the normalization $\nu_{(a,b)}:X_T\to X_c$. We considered above the map induced by pushforward. Now we take the pullback map 
\begin{equation}
    \label{eq:hatnu}
    \hat\nu_{(a,b)}:\Jac^{2g-2}(X_c)\to\Jac^{2g-2}(X_T)
\end{equation} on line bundles of degree $2g-2$. Under this map the generalized Jacobian $\Jac^{2g-2}(X_c)$ becomes a $(\CC^*)^{2g-2}$-bundle over $\Jac^{2g-2}(X_T)$. The morphism $\hat\nu_{(a,b)}$ does not extend to the compactification $\overline{\Jac}^{2g-2}(X_c)$ but each fixed fiber of $\hat\nu_{(a,b)}$ compactifies to $(\PP^1)^{2g-2}$ in $\overline{\Jac}^{2g-2}(X_c)$; cf. \cite{bhosle:1992}.

The family of spectral curves parameterized by $B^T_{ni}$ is composed by integral singular curves all of them normalized by $X_T$. In order to obtain a family of spectral curves equipped with the normalization, recall \eqref{eq X_c from b} and define
\[
\Cc \to \wt{B}^T
\]
to be the family of curves, whose fibre at $(a,b)$ is the spectral curve $\nu_{(a,b)}(X_T) = b(X_T) + a$ inside $|K|$, which amounts to $X_c$. By construction, $\Cc$ is equipped with the normalization
\[
\mathbf{n} = X_T \times \wt{B}^T \to \Cc,
\]
where $\mathbf{n}_{(a,b)} = \nu_{(a,b)}$. Note that $(a,b)$ and $(a,-b)$ in $\wt{B}^T$ define the same curve $X_c$ even though the normalization maps at these slices are not equal. Indeed, $\mathbf{n}_{(a,-b)} = \nu_{(a,-b)}$ is obtained from $\mathbf{n}_{(a,b)} = \nu_{(a,b)}$ after composing with the Galois involution. While we observe that $\Cc$ descends to a family of spectral curves parameterized by $B^T$, this is not the case of the relative morphism $\mathbf{n}$.

Set $\Cc_{ni}$ to be the restriction of $\Cc$ to $\wt{B}^T_{ni}$, the preimage of $B^T_{ni}$ in $\wt{B}^T$. Following \cite{altman&kleiman}, we can consider the corresponding relative Jacobian (of degree $2g-2$) $\mathcal{J}^{2g-2}\to \wt{B}^T_{ni}$ over $\wt{B}^T_{ni}$, which can be fiberwise compactified into the corresponding family $\overline{\mathcal{J}}^{2g-2}$ of compactified Jacobians. Let also $\mathcal{J}^{2g-2}_T=\Jac^{2g-2}(X_T)\times \wt{B}^T_{ni}\to \wt{B}^T_{ni}$ be the relative Jacobian for the constant family $X_T\times \wt{B}^T_{ni}$ of normalizations. Then we have a pullback map $\hat{\mathbf{n}}:\mathcal{J}^{2g-2}\to\mathcal{J}^{2g-2}_T$ induced by each $\hat\nu_{(a,b)}$ as in \eqref{eq:hatnu}. In addition, $\mathcal{J}^{2g-2}_T\to \wt B^T_{ni}$, being the trivial family, admits a section associated to each degree $g-1$ line bundle $M\in\Jac^{g-1}(X)$, given by 
\begin{equation}\label{eq:sections}
\sigma_M(a,b)=p^*M=\nu_{(a,b)}^*\pi_c^*M\in\Jac^{2g-2}(X_T).
\end{equation}

Note that under the spectral correspondence, one gets the map
\[
\jmath: \overline{\mathcal{J}}^{2g-2} \to \Mm
\]
which is $2:1$ as $\wt{B}^T_{ni} \to B^T_{ni}$ is $2:1$. The image of $\jmath$ is the subspace of $\Mm$ mapping to $B^T_{ni}$ under $h$. Taking preimage under $\hat{\mathbf{n}}$ of the image of $\sigma_M$ we get a subspace of $\mathcal{J}^{2g-2}$ (hence of $\Mm$) which maps onto $B^T_{ni}$. Define then $\Ss^M_T$ as the closure in $\mathcal{J}^{2g-2}$ and $B^T_{ni}$ of this subspace. In symbols, define
\begin{equation}\label{eq:defLag}
\Ss^M_T=\overline{\jmath \left (\hat{\mathbf{n}}^{-1}(\mathrm{Im}(\sigma_M)) \right )}\subset\Mm.
\end{equation}
Note that $h(\Ss^M_T)=B^T=h(\Mm^T)$, for any $M\in\Jac^{g-1}(X)$. Over every point $c=(2a,a^2-b^2) \in B^T_{ni}$, the restriction to the corresponding Hitchin fibre is
\[
\Ss^M_T\cap h^{-1}(c) = \overline{\hat{\nu}_{(a,b)}^{-1}(p^*M)} = \overline{\hat{\nu}_{(a,-b)}^{-1}(p^*M)},
\]
as $p^*M$ over $X_c$ is invariant under the Galois involution. When the distinction between $\nu_{(a,b)}$ and $\nu_{(a,-b)}$ is irrelevant, we shorten both to $\nu_c$, so that
\begin{equation}\label{eq:spectraldataLagrangian}
\Ss^M_T\cap h^{-1}(c)=\overline{\hat\nu_c^{-1}(p^*M)}\cong(\PP^1)^{2g-2}\subset\overline{\Jac}^{2g-2}(X_c)=h^{-1}(c).
\end{equation}

\begin{theorem}{\cite[Theorem 5.15]{FGOP}}\label{thm:Lagrangian}
The subvariety $\Ss^M_T$ of $\Mm$ is complex Lagrangian.
\end{theorem}

In \cite{FGOP}, $\Ss^M_T$ is defined in a different way, using Hecke transformations of Higgs bundles. Theorem \ref{thm:descriptionSs_T} below shows that both definitions yield the same subvariety $\Ss^M_T$. Before stating it, let us briefly recall what is a Hecke transformation along a simple divisor, which is the setting we are considering. 

Let $U$ be a rank $2$ vector bundle on $X$,  $p$ a point of $X$ and $W_p$ a $1$-dimensional vector subspace of the fiber $U_p$ of $U$ at $p$. Then, as a sheaf, $U_p/W_p$ is isomorphic to $\Oo_p$. Then the Hecke transform of $U$ at $p$ associated to $W_p$ is the kernel $E$ of the sheaf projection $U\to \Oo_p\cong U_p/W_p$. Now, if $D=p_1+\cdots+p_k$ is a simple effective divisor, after choosing a linear subspace $W_{p_i}\subset U_{p_i}$ for each $p_i$, then the Hecke transform of $U$ along $D$ associated to $\Vv'_D=(W_{p_1},\ldots,W_{p_k})$ is the kernel $E$ of the surjective sheaf map $U\to \Oo_D\cong \bigoplus_{i=1}^kU_{p_i}/W_{p_i}$. Then we have the exact sequence of sheaves
\[0\to E\to U\to \Oo_D\to 0,\]
and $E$ is locally free because it is torsion-free over the smooth curve $X$. Conversely, if $E$ is a rank $2$ subsheaf of $U$ such that $U/E$ is supported on a simple divisor, then it is a Hecke transformation of $U$ along the divisor supporting the torsion sheaf $U/E$. A Hecke transformation of Higgs bundles is a Hecke transformation of the underlying vector bundles which is compatible with the corresponding Higgs fields in the obvious way.

 This definition is equivalent to the one given above due to the next theorem, providing an alternative way to characterize Higgs bundles in the dense open subset $\Ss^M_T\cap h^{-1}(B^T_{ni})$ of $\Ss^M_T$.

\begin{theorem}{\cite[Theorem 5.5]{FGOP}}\label{thm:descriptionSs_T}
Let $c=(2a,a^2-b^2)\in B^T_{ni}$ with $a\in H^0(K)$ and $b\in H^0(KT)$. Recall that the effective simple divisor $\div(b)\in\Sym^{2g-2}(X)$ is the image under $\pi_c:X_c\to X$ of the $2g-2$ nodes of $X_c$. Let $(E,\varphi)\in \Mm$ and consider the rank $2$ and degree $2g-2$ Higgs bundle 
\begin{equation}\label{eq:Higgs-Heck-fixed-deg2g-2}
(M\oplus MT,\Phi),\ \ \text{with}\ \ \Phi=\begin{pmatrix} a & b \\ b & a \end{pmatrix}
\end{equation} where $M\in\Jac^{g-1}(X)$.
Then the following are equivalent:
\begin{enumerate}
\item $(E,\varphi)\in \Ss^M_T\cap h^{-1}(c)$;
\item $(E,\varphi)$ is a Hecke transformation of \eqref{eq:Higgs-Heck-fixed-deg2g-2} along the divisor $\div(b)$, that is, $E$ fits in an exact sequence $$0\to E\xrightarrow{\ \Psi\ } M\oplus MT \to \Oo_{\div(b)} \to 0$$ with $\varphi$ verifying $(\Psi\otimes 1_K)\varphi=\Phi\Psi$.
\end{enumerate}
Furthermore, a Higgs bundle $(E,\varphi)\in\Ss^M_T\cap h^{-1}(c)$ uniquely defines such a short exact sequence up to rescaling $\Psi$ and vice-versa.
\end{theorem}  

Notice that $KT\cong \Oo_X(\div(b))$ implies in particular that the underlying bundle of all Higgs bundles in $\Ss^M_T$ have fixed determinant isomorphic to $M^2K^{-1}$.

\

The main motivation for introducing the complex Lagrangian subvarieties $\Ss^M_T$ comes from the study of the behaviour under mirror symmetry of Narasimhan--Rammanan $\BBB$-branes. This study blows down to the following theorem, which constitutes the main result of \cite{FGOP}. 

\begin{theorem}{\cite{FGOP}}\label{thm:FM}
Over $B^T_{ni}$, the Fourier-Mukai transform of the Narasimhan-Ramanan $\BBB$-brane $\Delta^*\hat\Ll$, supported on $\Mm^T$, is a sheaf supported on the complex Lagrangian $\Ss^{M\otimes\Ll}_T$ of $\Mm$.
\end{theorem}

\section{Narasimhan--Ramanan branes and the wobbly locus} 
\label{sc NR and wobbly}

\subsection{Wobbly loci of $\BBB$-branes}
\label{sc BBB and wobbly}

In this section we study the intersection of the nilpotent cone $h^{-1}(0)$ with $\Mm^T$, the support of the $\BBB$-branes and  give a description of the wobbly loci therein.

First, observe that by Proposition 3.1 (ii) of \cite{NR}, any $(E,\varphi)\in \Mm^T$  satisfies that $E$ is semistable. So the following statement holds (recall that $\Nn$ denotes the moduli space of rank $2$ and degree $0$ vector bundles on $X$).

\begin{proposition}
If $(E,\varphi)\in\Mm^T$ then $\lim_{t\to 0}(E,t\varphi)=(E,0)\in\Nn$. 
\end{proposition}

Thus
\[
\Mm^T \cap h^{-1}(0) \cong \quotient{\Jac(X_T)}{\ZZ_2}\cong\Nn^T\subset\Nn,
\] 
where the line bundles $F \in \Jac(X_T)$ are sent to $E = p_* F$, where $p:X_T\to X$ is the étale $2$-cover defined by $T$, and where $\Nn^T$ is the analogous to $\Mm^T$ but in the action only on $\Nn$.

Let us now consider $\Ww^T:=\Mm^T\cap \Ww$, the intersection of $\Mm^T$ with the wobbly locus (necessarily of $\Nn$ by the previous proposition). It intersects the different components in subschemes denoted by
\[
\Ww_\delta^T:=\Mm^T\cap \Ww_\delta\subset \Ww^T:=\Mm^T\cap \Ww.
\]

In the remaining of the section we analyse $\Ww_\delta^T$. We recall from Section \ref{sc very stable} that $\delta \equiv 0 \! \mod 2$ and $0 \leq \delta \leq 2g-2$, and that $\Ww_\delta\subset\Ww_g$ for $\delta>g$. 

Denote by $\Nm$ the norm map on divisors associated to $p$. By definition, given a divisor $\sum a_ix_i$ in $X_T$ ($a_i\in\ZZ$, $x_i\in X_T$), $\Nm(\sum a_ix_i)$ is the divisor in $X$ defined by $\Nm(\sum a_ix_i)=\sum a_ip(x_i)$.

\begin{theorem}\label{thm:wobbly_BBB}
 Let $E\in\Nn$. Then $E\in\Ww_\delta^T$ if and only if $E=p_*\left(p^*F_0(R)\right)$ for some $F_0\in \Jac^{1-g+\delta/2}(X)$ and $R\in \Sym^{2g-2-\delta}(X_T)$ such that $R\cap\sigma^*R=\emptyset$ and moreover $KT(-\Nm(R))\cong\Oo(D)$ for some effective divisor $D\in\Sym^\delta(X)$. 
In particular, the subscheme $\Ww_\delta^T$ is non-empty for all possible values of $0\leq\delta\leq 2g-2$.   
\end{theorem}

\begin{proof}
Let $E=p_*F$ for some degree $0$ line bundle $F$ over $X_T$, so that $E\in\Nn^T$. Suppose $E$ is wobbly and let $\psi\in H^0(\End(E)\otimes K)\setminus\{0\}$ be nilpotent. Let $F_0\subset E$ be the line subbundle defined by $F_0=\ker(\psi)$ and $F_1=E/F_0$ be the quotient line bundle,
\[0\to F_0\to E\to F_1\to 0,\]
with $F_1\cong \det(E)F_0^{-1}$. Since $\psi$ factors through $\psi:F_1\longrightarrow F_0K$, we have that $F_1^{-1}F_0K\cong \Oo(D)$ for the effective divisor $D=\div(\psi)$ of even degree $\delta$ with $0\leq \delta\leq 2g-2$. So $E\in\Ww_\delta^T$.

Recall that $p^*E\cong F\oplus \sigma^*F$ where $\sigma:X_T\to X_T$ is the involution exchanging the sheets of the cover. So we have the exact sequence 
\begin{equation}\label{eq:ses}
0\to p^*F_0\to F\oplus \sigma^*F\to p^*F_1\to 0,    
\end{equation}
where $p^*F_1\cong F\otimes \sigma^*F\otimes p^*F_0^{-1}$.
Now, we have that either both the compositions $p^*F_0\hookrightarrow F\oplus \sigma^*F\twoheadrightarrow F$ and $p^*F_0\hookrightarrow F\oplus \sigma^*F\twoheadrightarrow \sigma^*F$ are non zero, or $p^*F_0\cong F\cong \sigma^*F$. In particular, 
\begin{equation}\label{eq:F0FR}
    F\cong p^*F_0\otimes\Oo_{X_T}(R)=p^*F_0(R)
\end{equation} for some effective divisor $R$ (possibly zero). From this it follows immediately that 
\[\deg(R)=-2\deg(F_0)=2g-2-\delta\] and that \[E\cong p_*(p^*F_0(R)).\] 
Now, it is well-known that $\det(p_*\Oo_{X_T}(R))\cong\Oo(\Nm(R))\otimes\det(p_*\Oo_{X_T})$. Using this and the fact that $p_*\Oo_{X_T}\cong\Oo\oplus T$, we get, by the projection formula, 
\begin{equation}\label{eq:detE}
    \det(E)\cong F_0^2\det(p_*\Oo_{X_T}(R))\cong F_0^2T(\Nm(R)).
\end{equation}
Since $F_0K(-D)\cong F_1\cong \det(E)F_0^{-1}$ we conclude from \eqref{eq:detE} that 
\[KT(-\Nm(R))\cong\Oo(D)\]
as claimed.

In addition, \eqref{eq:F0FR} also shows that the exact sequence \eqref{eq:ses} can be rewritten as
\[
0\to p^*F_0\xrightarrow{\ s\ } p^*F_0(R)\oplus p^*F_0(\sigma^*R)\to p^*F_1\to 0,
\] 
and since $p^*F_1$ is torsion-free, we conclude that the map $s$ cannot vanish simultaneously over $R$ and $\sigma^*R$. This implies that $R\cap\sigma^*R=\emptyset$.

Conversely, if $E=p_*\left(p^*F_0(R)\right)$ in the given conditions, then obviously $E\in\Nn^T$ and moreover $\det(E)\cong F_0^2T(\Nm(R))$ as in \eqref{eq:detE}. 

Now, take the diagonal embedding
\[
p^*F_0\hookrightarrow p^*E=p^*F_0(R)\oplus p^*F_0(\sigma^*R).
\]
It has torsion-free quotient because by assumption $R\cap\sigma^*R=\emptyset$. It is also an equivariant section, and so it descends to either $F_0\hookrightarrow E$ or $F_0T\hookrightarrow E$ (recall that $E\cong E\otimes T$), as the kernel of the pullback under $p$ is generated by $T$. Using that $\det(E)\cong F_0^2T(\Nm(R))$ and that by assumption $KT(-\Nm(R))\cong\Oo(D)$ for some $D\in\Sym^\delta(X)$, we conclude that the first case yields the short exact sequence 
\[0\to F_0\to E\to F_0K(-D)\to 0\]
while the second one produces
\[0\to F_0T\to E\to F_0KT(-D)\to 0.\] In both cases a non-zero section of $\Oo(D)\cong\Hom(F_0K(-D),F_0K)\cong\Hom(F_0KT(-D),F_0TK)$ gives rise to a non-zero nilpotent $\psi:E\to E\otimes K$ with kernel $F_0$ in the first case or $F_0T$ in the second one. In both cases the degree of this kernel is $1-g+\delta/2$ and so $E\in\Ww_\delta^T$. 

Finally, since the above construction works for any possible value of $\delta$, it follows that all the subschemes $\Ww_\delta^T$ are non-empty.
\end{proof}

\begin{remark}
When $R=0$, namely when $\Oo(D)\cong K$ or $\Oo(D)\cong KT$, we see that $E$ is wobbly if and only if it is strictly semistable. 
\end{remark}

\subsection{A stratification of the generic Hitchin fibres of the Lagrangian subvarieties}
\label{sc stratification}

In this section we construct a stratification of $\Ss_T^M\cap h^{-1}(c)$ for $c=(2a,a^2-b^2)\in B^T_{ni}$. As explained in Section \ref{sc C*flows of lagrangians} below, such stratification will allow us to compute the limit of the $\CC^*$-action of some points in $\Ss_T^M$, via the description of the maximal destabilizing bundle of the underlying vector bundle associated to these points.

Set $a\in H^0(K)$ and $b\in H^0(KT)$, such that $c=(2a,a^2-b^2)$ lies in $B^T_{ni}\subset H^0(K)\oplus H^0(K^2)$. By Theorem \ref{thm:descriptionSs_T}, any Higgs bundle in the intersection of the Hitchin fibre $h^{-1}(c)$ and our Lagrangian subvariety $\Ss^M_T$, $\Ee = (E,\varphi)$ in $\Ss^M_T \cap h^{-1}(c)$, is equipped with the short exact sequence
\begin{equation}\label{eq:Hecke}
0\to E\xrightarrow{\ \Psi\ } M\oplus MT \to \Oo_{\div(b)}\to 0
\end{equation} with $\varphi$ verifying $(\Psi\otimes 1_K)\varphi=\Phi\Psi$, and the short exact sequence \eqref{eq:Hecke} is uniquely defined up to scaling $\Psi$. Notice that if $p\in\div(b)$ then the $\Psi(E)_p$ is a $1$-dimensional linear subspace of $(M\oplus MT)_p$ and it coincides with the $1$-dimensional subspace $W_p$ defining the Hecke transformation at $p$ as explained before Theorem \ref{thm:descriptionSs_T}. Recall that $\Oo_X(\div(b))\cong KT$ and that moreover $\div(b)$ is reduced. Define
\begin{equation} \label{eq definition of D_Ee}
D_\Ee = \supp\Big(\quotient{M}{\Psi(E) \cap M}\Big).
\end{equation}
Then $D_\Ee$ is an effective subdivisor of $\div(b)$, hence it is reduced. It is also defined by considering the image of the restriction to $M$ of the map $M\oplus MT \to \Oo_{\div(b)}$. Indeed such restriction has image $\Oo_{D_\Ee}$. The points of $D_\Ee$ are precisely those $p\in\div(b)$ such that $\Psi(E)_p$ does not coincide with $M_p$. Hence  $D_\Ee=\div(b)$ for a generic Higgs bundle $\Ee\in\Ss^M_T\cap h^{-1}(c)$.

Given an effective subdivisor $0\leq D\leq \div(b)$, let us consider the locus of $\Ss_T^M\cap h^{-1}(c)$ given by those Higgs bundles yielding $D$ under \eqref{eq definition of D_Ee},
\[
\Vv(D) := \left \lbrace \Ee = (E,\varphi)\in\Ss_T^M \cap h^{-1}(c)  \textnormal{ such that } D_\Ee = D \right \rbrace.
\]
This naturally provides a stratification of the generic Hitchin fibres of our Lagrangian subvariety,
\begin{equation} \label{eq first stratification}
\Ss_T^M\cap h^{-1}(c)=\bigsqcup_{0\leq D\leq \div(b)}\Vv(D).
\end{equation}
To study the strata $\Vv(D)$ we need the following useful result.

\begin{lemma} \label{lm description of points in Vv}
Let $a\in H^0(K)$ and $b\in H^0(KT)$ with $c=(2a,a^2-b^2)\in B^T_{ni}\subset H^0(K)\oplus H^0(K^2)$. Pick $0 \leq D \leq \div(b)$ and consider a Higgs bundle $(E,\varphi)$ contained in $\Vv(D) \subset \Ss_T^M \cap h^{-1}(c)$. Then, $E$ fits in the commutative diagram
\begin{equation}\label{eq:commut-diag1}
\xymatrix{&&0\ar[d]&0\ar[d]&0\ar[d]&\\
          &0\ar[r]& M(-D)\ar[r]\ar[d]_{s}&E\ar[r]\ar[d]^-{\Psi}&MK^{-1}(D)\ar[r]\ar[d]^{b/s}&0\\
          &0\ar[r]&M\ar[r]\ar[d]&M\oplus MT\ar[r]\ar[d]&MT\ar[r]\ar[d]&0\\
          &0\ar[r]&\Oo_{D}\ar[d]\ar[r]&\Oo_{\div(b)}\ar[d]\ar[r]&\Oo_{\div(b)-D}\ar[r]\ar[d]&0\\
	  &&0&0&0,&}
\end{equation}
where $s\in H^0(\Oo_X(D))$, $b/s\in H^0(\Oo_X(\div(b)-D))$.
\end{lemma} 

\begin{proof}
Immediately after Theorem \ref{thm:descriptionSs_T} and the definition of $D_\Ee$ one has that $E$ fits into
\[
\xymatrix{&&0\ar[d]&0\ar[d]&&\\
          &0\ar[r]& M(-D_\Ee)\ar[r]\ar[d]_{s}&E\ar[d]^-{\Psi}&&\\
          &0\ar[r]&M\ar[r]\ar[d]&M\oplus MT\ar[r]\ar[d]&MT\ar[r]&0\\
          &0\ar[r]&\Oo_{D_\Ee}\ar[d]\ar[r]&\Oo_{\div(b)}\ar[d]&&\\
	  &&0&0.&&}
   \]
Recalling that $\det(E) \cong M^2K^{-1}$ one can easily complete the diagram to obtain \eqref{eq:commut-diag1} after identifying $D$ with $D_\Ee$, as one can do since $\Ee$ is a point of $\Vv(D)$.
\end{proof}

We set up some notation before stating the main result of this section. For a given effective divisor $D$ and section $s\in H^0(\Oo_X(D))$ with $\mathrm{div}(s)=D$, consider the double complex $C^\bullet_{D}:K(-2D)\xrightarrow{\,s\,}K(-D)$ given by multiplication by $s$. Associated to such complex, we have its first hypercohomology $\HH^1(C^\bullet_{D})$ vector space, defined as follows, using Dolbeault representatives. Consider the commutative square 
\begin{equation}\label{eq:commutsquare}
\xymatrix{&\Omega^0(K(-2D))\ar[d]_{\bar\partial}\ar[r]^{s}&\Omega^0(K(-D))\ar[d]^{\bar\partial}\\
          &\Omega^{0,1}(K(-2D))\ar[r]^{s}&\Omega^{0,1}(K(-D)).}
\end{equation}
Then, by definition,\[\HH^1(C^\bullet_{D})=\frac{\{(\sigma,\alpha)\in\Omega^{0,1}(K(-2D))\times \Omega^0(K(-D))\,|\, s\sigma=\bar\partial\alpha\}}{\{(\bar\partial x,sx)\in\Omega^{0,1}(K(-2D))\times \Omega^0(K(-D))\,|\, x\in\Omega^0(K(-2D))\}}\]
(since $X$ is a Riemann surface, $\sigma\in \Omega^{0,1}(K(-2D))$ is automatically $\bar\partial$-closed). Note that $\HH^1(C^\bullet_{D})$ is independent of the choice of the section $s$ and that if $D=0$ then $\HH^1(C^\bullet_{D})=\{0\}$.

\begin{proposition} \label{pr description of Vv}
Let $a\in H^0(K)$ and $b\in H^0(KT)$ are such that $c=(2a,a^2-b^2)\in B^T_{ni}\subset H^0(K)\oplus H^0(K^2)$. Take $\div(b)$ and fix an effective subdivisor $0\leq D\leq \div(b)$, then one has the identification
\[
\Vv(D) \cong \HH^1(C^\bullet_{D}), 
\]
so
\[
\dim \Vv(D) = \deg(D).
\]
\end{proposition}

\proof We will have to construct an element of $\HH^1(C^\bullet_{D})$ out of a given $(E,\varphi)\in\Vv(D)$ and conversely.

Let then $(E,\varphi)\in\Vv(D)$ so that $E$ is an extension of $MK^{-1}(D)$ by $M(-D)$.
 Let $\sigma\in H^1(K(-2D))$ be the corresponding extension class. Under the $C^\infty$ decomposition $M(-D)\oplus MK^{-1}(D)\cong E$, the holomorphic structure of $E$ is determined by the $\bar\partial$-operator $\bar\partial_E=\smtrx{ \bar\partial_1 & \sigma  \\ 0 & \bar\partial_2}$
with $\bar\partial_1$ (resp.~$\bar\partial_2$) the holomorphic structure of $M(-D)$ (resp.~$MK^{-1}(D)$) and where here $\sigma\in\Omega^{0,1}(\Hom(MK^{-1}(D),M(-D)))$ denotes a representative of the class $\sigma$.
We also have the holomorphic map $\Psi:E\to M\oplus MT$ determined by the Hecke transformation alone. Taking the above $C^\infty$-splitting of $E$, we write $\Psi=\smtrx{s & \alpha\\ 0 & b/s}$
for some $C^\infty$ section $\alpha\in\Omega^0(\Hom(MK^{-1}(D),M))=\Omega^0(K(-D))$. Let the holomorphic structure of $M\oplus MT$ being given by $\bar\partial_{M\oplus MT}=\left(\begin{smallmatrix} \bar\partial_3 &0 \\ 0 & \bar\partial_4
\end{smallmatrix}\right)$.
As $\Psi$ is holomorphic, $\bar\partial_{M\oplus MT}\Psi-\Psi\bar\partial_E=0$ which is equivalent to $
\bar\partial_3\alpha-\alpha\bar\partial_2=s\sigma$ and hence to $\bar\partial\alpha=s\sigma$ where $\bar\partial$ is the holomorphic structure on $K(-D)\cong\Hom(MK^{-1}(D),M)$ given on the right map of \eqref{eq:commutsquare}.

So out of $(E,\varphi)\in\Vv(D)$ we obtained the class in $\HH^1(C^\bullet_{D})$ represented by $(\sigma,\alpha)$,
\begin{equation}\label{eq:corresp1}
(E,\varphi)\in\Vv(D)\mapsto[(\sigma,\alpha)]\in\HH^1(C^\bullet_{D}).
\end{equation}
Note that the Higgs field $\varphi$ does not give any new information because it is completely determined by $\Psi$ and $\Phi$ via the condition $(\Psi\otimes 1_K)\varphi=\Phi\Psi$.

The converse inclusion is basically reversing the procedure just described, Suppose then that we are given an element $[(\sigma,\alpha)]\in \HH^1(C^\bullet_{D})$ represented by $(\sigma,\alpha)\in\Omega^{0,1}(K(-2D))\times\Omega^0(K(-D))$ such that
$s\sigma=\bar\partial\alpha$.
Take the cohomology class $\sigma\in H^1(K(-2D))=H^1(\Hom(MK^{-1}(D),M(-D)))$ and let $E$ be the associated extension of $MK^{-1}(D)$ by $M(-D)$,
so that its holomorphic structure is given by $\bar\partial_E=\smtrx{ \bar\partial_1 & \sigma  \\ 0 & \bar\partial_2}$. 
Define also $\Psi:E\to M\oplus MT$ as $\Psi=\smtrx{s & \alpha\\ 0 & b/s}$ in the $C^\infty$-splitting $M(-D)\oplus MK^{-1}(D)$ of $E$. Since both $s$ and $b/s$ are holomorphic, the condition $s\sigma=\bar\partial\alpha$ ensures that $\Psi$ is holomorphic. It is then clear that $E$ is a Hecke transformation of $M\oplus MT$ via $\Psi$ and fitting in \eqref{eq:commut-diag1}.

It remains to construct a Higgs field $\varphi:E\to E\otimes K$ on $E$. Recall that $M\oplus MT$ has the Higgs field $\Phi=\smtrx{a & b\\ b & a}$. Then $(\Psi\otimes 1_K)\varphi=\Phi\Psi$ shows that, in the $C^\infty$ decomposition $E\simeq M(-D)\oplus MK^{-1}(D)$, the Higgs field is given by
$\varphi=\smtrx{
-s\alpha+a & -\alpha^2+b^2/s^2\\ s^2 & s\alpha+a}$.
%
It is easy to see that  $\varphi$ is indeed a holomorphic map, that is, $\bar\partial_E\varphi-\varphi\bar\partial_E=0$.

From the previous description of $\varphi$, one can easily observe that $\tr(\varphi) = 2a$ and $\det(\varphi)= a^2 - b^2$. This shows that the spectral curve associated to $(E,\varphi)$ just constructed is the one associated to $c \in B^T_{ni}$, which is integral by assumption. This implies that $(E,\varphi)$ is stable, hence lies in $\Vv(D)$. 

Performing the same construction from a different choice $(\sigma+\bar\partial x,\alpha+sx)$ of the class $[(\sigma,\alpha)]\in\HH^1(C^\bullet_{D})$ yields a Higgs bundle $(E',\varphi')$ which is isomorphic to $(E,\varphi)$ by the isomorphism $\smtrx{1 & x\\ 0 & 1}:E'\to E$ in the $C^\infty$-decompositions $E'\simeq M(-D)\oplus MK^{-1}(D)\simeq E$.

So we have the correspondence 
\begin{equation}\label{eq:corresp2}
[(\sigma,\alpha)]\in\HH^1(C^\bullet_{D})\mapsto (E,\varphi)\in\Vv(D)
\end{equation}
and clearly \eqref{eq:corresp1} and \eqref{eq:corresp2} are inverse correspondences and prove the first part of the proposition.

The dimension claim is obvious from the isomorphism $\HH^1(C^\bullet_{D})\cong H^0(\Oo_{D})$, where $\Oo_{D}$ is the structure sheaf of the effective divisor $D\subset X$. The fact that these two spaces are isomorphic follows by taking the long exact sequence in cohomology and in hypercohomology associated to the complex $C^\bullet_{D}$.
\endproof


The previous proposition completes the description of the stratification introduced in \eqref{eq first stratification}.

\begin{corollary} \label{co stratification with Vv}
Let $a\in H^0(K)$ and $b\in H^0(KT)$ are such that $c=(2a,a^2-b^2)\in B^T_{ni}\subset H^0(K)\oplus H^0(K^2)$. Then $\Ss_T^M\cap h^{-1}(c)$ admits a stratification
\[\Ss_T^M\cap h^{-1}(c)=\bigsqcup_{0\leq D\leq \div(b)}\Vv(D)=\bigsqcup_{0\leq D\leq \div(b)}\HH^1(C^\bullet_{D}),\]
with $\Ss(\div(b))=\HH^1(C^\bullet_{\div(b)})$ the open stratum and $\Ss(0)=\{0\}$ the closed one.
\end{corollary}

Analogously to \eqref{eq definition of D_Ee}, given a Higgs bundle $\Ee = (E,\varphi) \in \Ss^M_T \cap h^{-1}(c)$, with $c \in B^T_{ni}$ as before, one can define the simple subdivisor of $\div(b)$,
\begin{equation} \label{eq definition of D'_Ee}
D'_\Ee =\supp\Big( \quotient{MT}{\Psi(E) \cap MT}\Big),
\end{equation}
having the same properties as $D_\Ee$ but now with respect to $MT$ instead of $M$.  Defining locus 
\[
\Vv'(D) := \left \lbrace \Ee = (E,\varphi)\in\Ss_T^M \cap h^{-1}(c)  \textnormal{ such that } D'_\Ee = D \right \rbrace.
\]
produces again a natural stratification,
\begin{equation} \label{eq second stratification}
\Ss_T^M\cap h^{-1}(c)=\bigsqcup_{0\leq D\leq \div(b)}\Vv'(D).
\end{equation}

The following Lemma \ref{lm description of points in Ww}, Proposition \ref{pr description of Ww} and Corollary \ref{co stratification with Ww} are the natural counterparts of  Lemma \ref{lm description of points in Vv}, Proposition \ref{pr description of Vv} and Corollary \ref{co stratification with Vv}, respectively, and their proofs are exactly the same.

\begin{lemma} \label{lm description of points in Ww}
Let $a\in H^0(K)$ and $b\in H^0(KT)$ with $c=(2a,a^2-b^2)\in B^T_{ni}\subset H^0(K)\oplus H^0(K^2)$. Pick $0 \leq D \leq \div(b)$ and consider a Higgs bundle $(E,\varphi)$ contained in $\Vv'(D) \subset \Ss_T^M \cap h^{-1}(c)$. Then, $E$ fits in the commutative diagram
\begin{equation}\label{eq:commut-diag2}
\xymatrix{&&0\ar[d]&0\ar[d]&0\ar[d]&\\
          &0\ar[r]& MT(-D)\ar[r]\ar[d]_{s}&E\ar[r]\ar[d]^-{\Psi}&MTK^{-1}(D)\ar[r]\ar[d]^{b/s}&0\\
          &0\ar[r]&MT\ar[r]\ar[d]&M\oplus MT\ar[r]\ar[d]&M\ar[r]\ar[d]&0\\
          &0\ar[r]&\Oo_{D}\ar[d]\ar[r]&\Oo_{\div(b)}\ar[d]\ar[r]&\Oo_{\div(b)-D}\ar[r]\ar[d]&0\\
	  &&0&0&0,&}
\end{equation}
where $s\in H^0(\Oo_X(D))$, $b/s\in H^0(\Oo_X(\div(b)-D))$.
\end{lemma} 


\begin{proposition} \label{pr description of Ww}
Let $a\in H^0(K)$ and $b\in H^0(KT)$ are such that $c=(2a,a^2-b^2)\in B^T_{ni}\subset H^0(K)\oplus H^0(K^2)$. Take $\div(b)$ and fix an effective subdivisor $0\leq D\leq \div(b)$, then one has the identification
\[
\Vv'(D) \cong \HH^1(C^\bullet_{D}), 
\]
so
\[
\dim \Vv'(D) = \deg(D).
\]
\end{proposition}


\begin{corollary} \label{co stratification with Ww}
Let $a\in H^0(K)$ and $b\in H^0(KT)$ are such that $c=(2a,a^2-b^2)\in B^T_{ni}\subset H^0(K)\oplus H^0(K^2)$. Then $\Ss_T^M\cap h^{-1}(c)$ admits a stratification
\[
\Ss_T^M\cap h^{-1}(c)=\bigsqcup_{0\leq D\leq \div(b)}\Vv'(D)=\bigsqcup_{0\leq D\leq \div(b)}\HH^1(C^\bullet_{D}),
\]
with $\Ss(\div(b))=\HH^1(C^\bullet_{\div(b)})$ the open stratum and $\Ss(0)=\{0\}$ the closed one.
\end{corollary}

\begin{remark}
Even if $\Vv(D)$ and $\Vv'(D)$, for a fixed effective divisor $D \leq \div(b)$, are both isomorphic to $\HH^1(C^\bullet_{D})$, they parameterize different subsets of $\Ss_T^M\cap h^{-1}(c)$. For instance, 
\begin{equation}\label{eq:Hitsec1}
\Vv(0) = \left(M\oplus MK^{-1},\left(\begin{smallmatrix}a & b^2 \\ 1 & a\end{smallmatrix}\right)\right),
\end{equation}
whereas
\begin{equation}\label{eq:Hitsec2}
\Vv'(0)=\left(MT\oplus MTK^{-1},\left(\begin{smallmatrix}a & b^2 \\ 1 & a\end{smallmatrix}\right)\right).
\end{equation}
\end{remark}

It follows by all of the above that one can obtain a third statification by intersecting the previous two,
\[
\Ss_T^M\cap h^{-1}(c)=\bigsqcup_{0\leq D,D' \leq \div(b)}\Vv(D) \cap \Vv'(D').
\]


The next statement will allow us to analyize when we get an empty intersection of the strata $\Vv(D) \cap \Vv'(D')$.

\begin{lemma} \label{lm comparison two stratifications}
Let $a\in H^0(K)$ and $b\in H^0(KT)$ with $c=(2a,a^2-b^2)\in B^T_{ni}$. Consider two effective divisors $0 \leq D, D' \leq \div(b)$. Then 
\[
\Vv(D) \cap \Vv'(D') = \emptyset
\]
unless 
\[
\div(b) \leq D + D'. 
\]
In particular, $\Vv(D) \cap \Vv'(D') = \emptyset$ when
\[
0 \leq \deg(D), \deg(D') < g-1.
\]
\end{lemma}

\begin{proof}
Take any Higgs bundle $\Ee = (E,\varphi) \in \Ss_T^M\cap h^{-1}(c)$. Observing \eqref{eq definition of D_Ee}, we note that the points parameterized by $D_\Ee$ correspond to those points $x \in \div(b)$ such that $\Psi(E)|_x$ trivially intersects $M|_x$. Conversely, the points in $\div(b) - D_\Ee$ are those points $x' \in \div(b)$ where $\Psi(E)|_{x'} = M|_{x'}$. An analogous analysis can be applied to $D'_\Ee$ and $y \in D'_\Ee$ if and only if $\Psi(E)|_y \cap MT|_y = 0$ and $y' \in \div(b) - D'_\Ee$ if and only if $\Psi(E)|_{y'} = MT|_{y'}$. From this discussion, one observes that 
\[
(\div(b) - D_\Ee) \cap (\div(b) - D'_\Ee) = \emptyset.
\]
Hence, the sum of $\div(b) - D_\Ee$ and $\div(b) - D'_\Ee$ is an effective subdivisor of $\div(b)$, 
\[
\left( \div(b) - D_\Ee \right ) + \left ( \div(b) - D'_\Ee \right ) \leq \div(b). 
\]
This yields
\[
\div(b) \leq D_\Ee + D'_\Ee,
\]
and the proof of the lemma follows naturally from this statement.
\end{proof}

\begin{remark}
Note that $M|_x$ is subspace of dimension $1$ in the ambience vector space $(M \oplus MT)|_x$ of dimension $2$. Since for every $x \in \div(b)$ one has that $\Psi(E)|_x$ has dimension $1$, it follows that $\Psi(E)|_x \cap M|_x$ is generically trivial as so is the generic intersection of two subspaces of dimension $1$ in a vector space of dimension $2$. A similar analysis can be performed on the points of $D'_\Ee$. Indeed, the generic points are those parameterized by $\Vv(\div(b)) \cap \Vv'(\div(b))$, which is dense in $\Ss_T^M\cap h^{-1}(c)$.
\end{remark}


\subsection{Limits of $\CC^*$-flows along the Lagrangians}
\label{sc C*flows of lagrangians}

Our aim in this section is to study the Lagrangian subvarieties $\Ss^M_T$ at the nilpotent cone. Since $\Ss^M_T$ is constructed as the closure of a subvariety defined over the locus $B^T_{ni}$, we study the limits of the $\CC^*$-action on the points $\Ss^M_T$ inside the nilpotent cone. Our analysis shows that the Lagrangians $\Ss^M_T$ hit the wobbly locus.

\begin{proposition}
Let $M\in\Jac^{g-1}(X)$ and $\Ss^M_T$ be the Lagrangian of $\Mm$ defined in \eqref{eq:defLag}. Then, $\Ss^M_T\cap h^{-1}(B^T_{ni})$ is $\CC^*$-invariant. 
\end{proposition}

\proof
Take a generic Higgs bundle $(E,\varphi)$ in $\Ss^M_T$ in the sense that $(E,\varphi)\in\Ss^M_T\cap h^{-1}(c)$ for some $c=(2a,a^2-b^2)\in B^T_{ni}$ with $a\in H^0(K)$ and $b\in H^0(KT)$. Then by Theorem \ref{thm:descriptionSs_T}, it is a Hecke transformation of $(M\oplus MT,\Phi)$ given in \eqref{eq:Higgs-Heck-fixed-deg2g-2}. For $t\in\CC^*$, define $tc=(2ta,(ta)^2-(tb)^2)$ and note that $tc\in B^T_{ni}$ because $\Mm^T$ is $\CC^*$-invariant and also because $X_{tc}$ is still nodal and integral. Moreover, $(\Psi\otimes 1_K)t\varphi=t\Phi\Psi$ so $(E,t\varphi)\in h^{-1}(tc)$ is a Hecke transformation of $(M\oplus MT,t\Phi)$ along $D_{tb}=\div(b)$, hence again Theorem \ref{thm:descriptionSs_T} tells us that $(E,t\varphi)\in\Ss^M_T\cap h^{-1}(tc)$.
\endproof


Let us now study the fixed point loci of the $\CC^*$-action on the Lagrangian subvarieties $\Ss^M_T$. Recall that, as we have seen in Section \ref{sc very stable}, this amounts to study the maximal destabilizing subbundle of the underlying vector bundle of those Higgs bundles parameterized by $\Ss^M_T$.

\begin{proposition} \label{pr stability of the underlying bundle}
Consider $\Ee = (E,\varphi)\in\Ss^M_T\cap h^{-1}(c)$, with $c=(2a,a^2-b^2) \in B^T_{ni}$ for a certain $a \in H^0(K)$ and $b \in H^0(KT)$. If $E$ is not semistable, if and only if one of the two (mutually exclusive) conditions holds:
\begin{enumerate}[(a)]
    \item \label{it D} $\deg(D_\Ee) < g-1$ and the maximal destabilizing subbundle of $E$ is $M(-D_\Ee)$ and the composition $\phi:M(-D_\Ee)\to E\xrightarrow{\ \varphi\ }E\otimes K\to M(D_\Ee)$ is $s_\Ee^2$, where $s_\Ee$ is a section of $\Oo_X(D_\Ee)$ that vanishes in $D_\Ee$; 

\noindent
or
    
    \item \label{it D'} $\deg(D'_\Ee) < g-1$ and the maximal destabilizing subbundle of $E$ is $LMT(-D'_\Ee)$ and the composition $\phi:MT(-D'_\Ee)\to E\xrightarrow{\ \varphi\ }E\otimes K\to MT(D'_\Ee)$ is $(s'_\Ee)^2$, where $s'_\Ee$ is a section of $\Oo_X(D'_\Ee)$ that vanishes in $D'_\Ee$.
\end{enumerate}

Furthermore, $E$ is strictly semistable if and only if $\deg(D_\Ee) = g-1$ (and $M(-D_\Ee) \subset E$) or $\deg(D'_\Ee) = g-1$ (and $MT(-D'_\Ee) \subset E$). 
\end{proposition}

\begin{proof}
First observe that \eqref{it D} and \eqref{it D'} are mutually exclusive due to Lemma \ref{lm comparison two stratifications}.

If $\deg(D_\Ee) < g-1$ ({\it resp.} $\deg(D_\Ee) < g-1$) one immediately sees after Lemma \ref{lm description of points in Vv} ({\it resp.} Lemma \ref{lm description of points in Ww}) that $M(-D_\Ee)$ ({\it resp.} $MT(-D'_\Ee)$) is a destabilizing subbundle, so $E$ is not semistable. Conversely, suppose that $E$ is not semistable with maximal destabilizing subbundle $L$. Since the determinant of $E$ is $M^2K^{-1}$, its underlying vector bundle is equipped with the short exact sequence,
\[
0 \to L \to E \to M^2K^{-1}L^{-1} \to 0,
\]
and the composition $\phi:L\to E\xrightarrow{\ \varphi\ }E\otimes K\to M^2L^{-1}$ is non-zero, as otherwise it would contradict the semistability of $\Ee$. 

Using again Theorem \ref{thm:descriptionSs_T} to see $(E,\varphi)$ as a Hecke transformation of $(M\oplus MT,\Phi)$, it follows that $L$ is a subsheaf of $M\oplus MT$ with quotient having torsion over a subdivisor $D_L\leq \div(b)$. Saturating $L$ in $M\oplus MT$ produces the line subbundle $L(D_L)\subset M\oplus MT$. From this, and if $s_L$ is a section with divisor $D_L$, we get the following diagram, analogue to \eqref{eq:commut-diag1} and \eqref{eq:commut-diag2},
\begin{equation}\label{eq:commut-diag3}
\xymatrix{&&0\ar[d]&0\ar[d]&0\ar[d]&\\
          &0\ar[r]& L\ar[r]\ar[d]_{s_L}&E\ar[r]\ar[d]^-{\Psi}&M^2K^{-1}L^{-1}\ar[r]\ar[d]^{b/s_L}&0\\
          &0\ar[r]&L(D_L)\ar[r]\ar[d]&M\oplus MT\ar[r]\ar[d]&M^2TL^{-1}(-D_L)\ar[r]\ar[d]&0\\
          &0\ar[r]&\Oo_{D_L}\ar[d]\ar[r]&\Oo_{\div(b)}\ar[d]\ar[r]&\Oo_{\div(b)-D_L}\ar[r]\ar[d]&0\\
	  &&0&0&0,&}
\end{equation} where again $D_L$ records precisely the points $p\in \div(b)$ where $L(D_L)|_p$ does not coincide with $\Psi(E)|_p$ as a subsheaf of $(M \oplus MT)|_p$.

Let $\pr:(M\oplus MT)\otimes K\to M^2TL^{-1}K(-D_L)$ be the projection induced by the one on the second line above. Since $(\Psi\otimes 1_K)\varphi=\Phi\Psi$, it follows that $\pr\circ(\Psi\otimes 1_K)\circ\varphi|_L=\pr\circ\Phi\circ\Psi|_L$ and this equivalent to $(b/s_L)\phi=bs_L$ (where here $b/s_L:M^2L^{-1}\to M^2TL^{-1}K(-D_L)$). So $\phi=s_L^2$ and moreover 
\begin{equation} \label{eq L up to square root}
\Oo(2D_L)\cong M^{2}L^{-2},
\end{equation}
so
\[
L \cong M(-D_L) T'
\]
for some $T' \in \Jac(X)[2]$. This implies that $L(D_L) \cong MT'$ is a subbundle of $M \oplus MT$, or equivalently, $T'$ is a subbundle of $\Oo_X \oplus T$ which is only the case if $T' \cong \Oo_X$ or $T' \cong T$. In the first case, one has that 
\[
L \cong M(-D_L)
\]
and 
\[
D_L = D_\Ee,
\]
while in the second case 
\[
L \cong MT(-D_L)
\]
and 
\[
D_L = D'_\Ee.
\]
The first statement follows easily from the previous identifications. 

Finally, observe that when $\deg(D_\Ee) = g-1$, Lemma \ref{lm description of points in Vv} implies that $E$ is strictly semistable with $M(-D_\Ee) \subset E$ having trivial degree (and slope). Analogously, when $\deg(D'_\Ee) = g-1$, Lemma \ref{lm description of points in Ww} implies that $E$ is strictly semistable with $MT(-D'_\Ee) \subset E$ having trivial degree (and slope). Conversely, if $E$ is strictly semistable with $L \subset E$ topologically trivial, it follows from the previous analysis that $L(D_L)$ is either $M$ or $MT$. In the first case, this implies that $D_L = D_\Ee$ and $M(-D_\Ee) \subset E$ while in the second case, one has $D_L = D'_\Ee$ and $MT(-D'_\Ee) \subset E$. This proves the second claim and concludes the proof. 
\end{proof}

Recall from Section \ref{sc stratification} the stratifications $\bigsqcup_{D \leq \div(b)} \Vv(D)$ and $\bigsqcup_{D \leq \div(b)} \Vv'(D)$ of $\Ss^M_T \cap h^{-1}(c)$, the restriction of our Lagrangian subvarieties to a generic Hitchin fibre. Consider a Higgs bundle $\Ee = (E, \varphi)$ in the previous intersection. Thanks to Proposition \ref{pr stability of the underlying bundle} one has that the underlying vector bundle $E$ of $\Ee$ is unstable if and only if it lies in some $\Vv(D)$ or $\Vv'(D)$ for $D \leq \div(b)$ with $\deg(D) < g-1$. In view of this, we consider the set of subdivisors of $\div(b)$ with length smaller than $g-1$, 
\[
U_b = \left \lbrace   D \leq \div(b) \textnormal{ such that } \deg(D) < g-1  \right \rbrace.
\]
Hence, $E$ is semistable provided it is not contained in the union of $\bigsqcup_{D \in U_b} \Vv(D)$ and $\bigsqcup_{D \in U_b} \Vv'(D)$. Also, observe that it follows from Lemma \ref{lm comparison two stratifications} that the previous sets do not intersect. Furthermore, if $\Ee \in \Vv(D)$ or $\Ee \in \Vv'(D)$ for any $D \in U_b$, we immediately knows that its maximal destabilizing subbundle is $M(-D)$ or $MT(-D)$, respectively. Hence, knowing in which stratum lies $\Ee$ we know the maximal destabilizing bundle of its underlying vector bundle $E$, hence we know the limit of $\Ee$ under the $\CC^*$-action. Let us summarize all this in the following result.

\begin{corollary} \label{co C* limit}
Consider $\Ee = (E,\varphi)\in\Ss^M_T\cap h^{-1}(c)$, with $c=(2a,a^2-b^2) \in B^T_{ni}$ for some $a \in H^0(K)$ and $b \in H^0(KT)$. Then, 
\begin{itemize}
     \item if $\Ee = (E,\varphi)$ lies in the complement of $\bigsqcup_{D \in U_b} \Vv(D) \sqcup \bigsqcup_{D' \in U_b} \Vv'(D')$, then $E$ is semistable and 
     \[
     \lim_{t \to 0} (E, t \varphi) = (E,0);
     \]

    \item  if $\Ee = (E,\varphi)$ lies in $\Vv(D)$ ({\it resp.} in $\Vv'(D)$) for some $D \in U_b$, then $E$ is unstable and 
     \begin{eqnarray}\label{eq:limit_1}
    \lim_{t\to 0}(E,t\varphi)= \left(M(-D)\oplus MK^{-1}(D),\smtrx{0 & 0\\ s_D^2 & 0}\right), 
    \\ \nonumber
    \left ( {\it resp.} \quad \lim_{t\to 0}(E,t\varphi)=\left(MT(-D)\oplus MTK^{-1}(D),\smtrx{0 & 0\\ s_D^2 & 0}\right) \right ),
    \end{eqnarray}
    where $s_D$ is a non-zero section of $\Oo_X(D)$ with divisor $D$.
    
\end{itemize}
\end{corollary}

We can now address the main result of this section.

\begin{theorem}\label{thm:wobblylimits}
Consider $\Ee = (E,\varphi)\in\Ss^M_T\cap h^{-1}(c)$, with $c=(2a,a^2-b^2) \in B^T_{ni}$ for a certain $a \in H^0(K)$ and $b \in H^0(KT)$. Then, 
\begin{enumerate}
\item \label{it semistable wobbly} if $E$ is semistable ({\it i.e.} $\Ee = (E,\varphi)$ lies in the complement of $\bigsqcup_{D \in U_b} \Vv(D) \sqcup \bigsqcup_{D' \in U_b} \Vv'(D')$), then $\lim_{t\to 0}(E,t\varphi)$ is wobbly, that is, $E\in\Nn$ is wobbly.

\item \label{it unstable wobbly} If $E$ is unstable with $M(-D)\subset E$ maximal destabilizing subbundle for $0 < D \in U_b$ ({\it i.e.} $\Ee$ lies in $\Vv(D)$ or $\Vv'(D)$ for a non-zero $D \in U_b$), then $\lim_{t\to 0}(E,t\varphi)$ is wobbly, that is, \eqref{eq:limit_1} is wobbly.

\item \label{it very stable} If $E$ is unstable with $M\subset E$ maximal destabilizing subbundle $M$ ({\it i.e.} $\Ee$ lies in $\Vv(0)$ or $\Vv'(0)$ for a non-zero $D$), then $\lim_{t\to 0}(E,t\varphi)$ is very stable. 
\end{enumerate}  
\end{theorem}

\begin{proof}
Suppose $E$ is semistable. Observe that $E$ is obtained via an extension
\[
0 \to M(-D) \to E \to MK^{-1}(D) \to 0,
\]
(or with $M$ replaced by $MT$) for some $0\leq D\leq\div(b)$. Hence to prove that $E$ is wobbly, it is enough to show that there is a non-zero holomorphic map $MK^{-1}(D)\to M(-D)K$, because then the composition $E\to MK^{-1}(D)\to M(-D)K\to E\otimes K$ will be a non-zero nilpotent Higgs field on $E$. But $H^0(\Hom(MK^{-1}(D),M(-D)K))=H^0(\Oo_X(2\div(b)-2D))\neq 0$ because $\div(b)-D$ is effective. This proves \eqref{it semistable wobbly}. 

Item \eqref{it unstable wobbly} follows from Corollary \ref{co C* limit} and the criterion \cite[Theorem 1.2]{HH}, noting that $s_D^2$ has double zeros (because $D>0$) and that the Higgs bundles of the form \eqref{eq:limit_1} are stable. 

Regarding the last item, we have that \eqref{eq:limit_1} equal
\begin{eqnarray*}
    \lim_{t\to 0}(E,t\varphi)= \left(M\oplus MK^{-1},\smtrx{0 & 0\\ 1 & 0}\right), 
    \\ \nonumber
    \left ( {\it resp.} \quad \lim_{t\to 0}(E,t\varphi)=\left(MT\oplus MTK^{-1},\smtrx{0 & 0\\ 1 & 0}\right) \right ),
    \end{eqnarray*}
    so they are fixed points with destabilzing subbundle of $E$ of maximal degree, hence they are very stable by Remark \ref{rmk:nomaximalwobbly}, hence proving \eqref{it very stable}. Indeed, in such case, $\Ee$ must be of the form \eqref{eq:Hitsec1} (resp. \eqref{eq:Hitsec2}) i.e.~must lie in a Hitchin section (associated to $M$) of $h$.
\end{proof}

\begin{remark}\mbox{}
\begin{enumerate}
\item According to \cite[Theorem 1.2]{HH} a rank $2$ Higgs bundle $\left(L_1\oplus L_2,\begin{psmallmatrix}0 & 0\\ \phi & 0\end{psmallmatrix}\right)$ with $\phi\neq 0$ is wobbly if and only if $\phi$ has a multiple zero. Hence one concludes that the wobbly Higgs bundles (with non-vanishing Higgs field) arising as $\CC^*$-limits along the Lagrangians $\Ss_T^M$ are special in the sense that \emph{all} zeros are multiple.
\item Theorem \ref{thm:wobblylimits} says that the intersection of the Lagrangians $\Ss^M_T$ with the nilpotent cone i.e.~the closure of $\Ss^M_T\cap h^{-1}(B^T_{ni})$ in $h^{-1}(0)$ is not just given by the limit point under the $\CC^*$-flows in $\Ss^M_T\cap h^{-1}(B^T_{ni})$. Indeed, since they are wobbly, the upward flow of such limit points contains $\CC^*$-orbits inside the nilpotent cones. Then those orbits in $h^{-1}(0)$ whose limits (to $0$ and $\infty$) are fixed points which are limits of Higgs bundles in $\Ss^M_T\cap h^{-1}(B^T_{ni})$ belong to $\Ss_T^M$. It would be interesting to have an explicit description of these nilpotent $\CC^*$-orbits in $\Ss^M_T$. Indeed, this would give a complete description of $\Ss^M_T\cap h^{-1}(0)$. We describe these flows in work in progress.
\item Upon varying $c=(2a,a^2-b^2)$ in $B^T_{ni}$, \eqref{eq:Hitsec1} and \eqref{eq:Hitsec2} give rise to two different sections of $h|_{h^{-1}(B^T_{ni})}$ which extend to $B^T$, so take values in $\Ss^M_T$. In particular, choosing a spin structure $K^{1/2}$, we see that the Lagrangian $S_T^{K^{1/2}}$ contains (the restrictions to $B^T$ of) two classical Hitchin sections, namely $\left(K^{1/2}\oplus K^{-1/2},\left(\begin{smallmatrix}a & b^2 \\ 1 & a\end{smallmatrix}\right)\right)$ and $\left(K^{1/2}T\oplus K^{-1/2}T,\left(\begin{smallmatrix}a & b^2 \\ 1 & a\end{smallmatrix}\right)\right)$, associated to the two spin structures $K^{1/2}$ and $K^{1/2}T$. By taking the locus in $B^T$ with $a=0$ these two sections lie also in the moduli space of $\SL(2,\CC)$-Higgs bundles. In such moduli space, the $\CC^*$-fixed point locus of the form \eqref{eq:nonzerofixedpoint} with $\deg(L_1)=g-1$ is a finite set because $L_1$ must be a square root of $K$. Hence it has $2^{2g}$ points.  From above we see that $S_T^{K^{1/2}}$ intersects exactly two of those $2^{2g}$ points: $\left(K^{1/2}\oplus K^{-1/2},\left(\begin{smallmatrix}0 & 0 \\ 1 & 0\end{smallmatrix}\right)\right)$ and $\left(K^{1/2}T\oplus K^{-1/2}T,\left(\begin{smallmatrix}0 & 0 \\ 1 & 0\end{smallmatrix}\right)\right)$.

\end{enumerate}
\end{remark}


\begin{thebibliography}{HMDP}

\bibitem[AK]{altman&kleiman}
A. Altman, S. Kleiman, {\it Compactifying the Picard scheme. II},  {Amer. J. Math.} \textbf{101} (1979), no. 1, 10--41.

\bibitem[Ar]{arinkin}
D. Arinkin, {\it Autoduality of compactified Jacobians for curves with plane singularities}, {J. Algebraic Geom.} \textbf{22} (2013), 363--388.

\bibitem[BNR]{BNR}
A. Beauville, M. S. Narasimhan, S. Ramanan, {\it Spectral curves and the generalized theta divisor}, {J. Reigne Angew. Math} \textbf{398} (1989), 169--179.

\bibitem[Bh]{bhosle:1992}
{U. Bhosle}, Generalised parabolic bundles and applications to torsionfree sheaves on nodal curves, \emph{Arkiv f\"ur Mat.} (2) \textbf{30} (1992), 187--215.


\bibitem[C2]{cook:1998}
{P. R. Cook}, {\it Compactified Jacobians and curves with simple singularities},  Algebraic Geometry
(Catania, 1993/Barcelona, 1994), 37--47, Lecture Notes in Pure and Appl. Math., \textbf{200}, Marcel Dekker, 1998.

\bibitem[DP]{DP} 
{R. Donagi, T. Pantev}, {\it Langlands duality for Hitchin systems},
{Invent. Math.} \textbf{189}, no. 3 (2012), 653--735.


\bibitem[FGOP]{FGOP}
E. Franco, P. Gothen, A. Oliveira, A. Peón-Nieto, {\it Unramified covers and branes on the Hitchin system}, {Adv. Math.} \textbf{377} (2021), 107493.


\bibitem[FP]{FP}
E. Franco, A. Peón-Nieto, {\it Branes on the singular locus of the Hitchin system via Borel and other parabolic subgroups}, {Math. Nach.} \textbf{296} (2023), 1803--1841.



\bibitem[GR1]{garciaprada-ramanan-2}
O. García-Prada, S. Ramanan, {\it Involutions of rank 2 Higgs bundle moduli spaces}, Geometry and Physics: Volume II, A Festschrift in honour of Nigel Hitchin, Editors: J.E. Andersen, A. Dancer and O. Garcia-Prada, Oxford University Press, 2018.

\bibitem[GR2]{garciaprada-ramanan}
O. García-Prada, S. Ramanan, {\it Involutions and higher order automorphisms of Higgs moduli spaces}, Proc. London Math. Soc. \textbf{119} (2019), 681--732.

\bibitem[GO1]{gothen-oliveira:2013}
P. Gothen, A. Oliveira,  {\it The singular fiber of the Hitchin map}, {Int. Math. Res. Not.}, \textbf{2013} No. 5 (2013), 1079--1121.

\bibitem[HT]{HT} 
T. Hausel, M. Thaddeus, {\it Mirror symmetry, Langlands duality, and the Hitchin system}, {Invent. Math.} \textbf{153} (2003), 197--229.

\bibitem[HH]{HH} 
T. Hausel, N. Hitchin, {\it Very stable Higgs bundles, equivariant multiplicity and mirror symmetry}, {Invent. Math.} \textbf{228} (2022), 893--989.

\bibitem[Hi1]{hitchin-self} 
N. J. Hitchin, {\it The self-duality equations on a Riemann
surface}, {Proc. London Math. Soc.} (3), \textbf{55}, 1, (1987), 59--126.


\bibitem[Hi2]{hitchin_duke} 
N. J. Hitchin, {\it Stable bundles and integrable systems}, {Duke Math. J.} \textbf{54}, Number 1 (1987), 91--114. 


\bibitem[Hi3]{hitchin_G2} 
N. J. Hitchin, {\it Langlands Duality And $G_2$ Spectral Curves}, {Quart. J. Math.} \textbf{58}, 3, (2007), 319--344. 

\bibitem[KW]{kapustin&witten}
A. Kapustin, E. Witten, {\it Electric-magnetic duality and the geometric Langlands program}, {Comm. Number Theory Phys.} \textbf{1} (2007), 1--236.    


\bibitem[Ko]{kont:1995} 
M. Kontsevich, {\it Homological Algebra of Mirror Symmetry}, In: Chatterji, S.D. (eds) Proceedings of the International Congress of Mathematicians. Birkhäuser, Basel, 1995.

\bibitem[Na]{nasser:2005} F. Nasser, {\it Torsion Subgroups of Jacobians Acting on Moduli Spaces of Vector Bundles}, {PhD Thesis}, University of Aarhus, 2005.

\bibitem[NR]{NR}
M. S. Narasimhan, S. Ramanan, {\it Generalized Prym varieties as fixed points}, {J. Indian Math. Soc.} \textbf{39} (1975), 1--19.

\bibitem[PaP]{PalPauly}
S. Pal, C. Pauly, {\it The wobbly divisors of the moduli space of rank-2 vector bundles}, Advances in Geometry 21 (4), 473--482.

\bibitem[PPe]{pauly&peon-nieto}
C. Pauly, A. Peón-Nieto, {\it Very stable bundles and properness of the Hitchin map.}, {Geom. Dedicata} \textbf{198} (2019) 143--145.

\bibitem[Si1]{simpson-cstar}  C.T. Simpson, {\it Higgs bundles and local systems},
Pub. Math. IHES, Tome 75 (1992), 5--95.

\bibitem[Si2]{Si94}
C.T. Simpson, {\it Moduli of representations of the fundamental group of a smooth projective variety {I}}, {Publ. Math. I.H.E.S.} \textbf{79} (1994), 47--129.

 
\bibitem[SYZ]{strominger-yau-zaslow:1993}
{A. Strominger, S.-T. Yau, E. Zaslow}, {\it Mirror symmetry is T-duality}, 
{Nuclear Phys. B} \textbf{479} (1996), 243--259.

\end{thebibliography}
\end{document}